\documentstyle[12pt]{article}
     \def\proof{\vskip 3mm \noindent{\bf Proof:}\hskip10pt}
     \def\QED{\hfill $\Box$ \smallskip}

     \font\tenmath=msbm10 scaled 1200
     \font\sevenmath=msbm7 scaled 1200
     \font\fivemath=msbm5 scaled 1200
     \newfam\mathfam \textfont\mathfam=\tenmath
     \scriptfont\mathfam=\sevenmath \scriptscriptfont\mathfam=\fivemath
     \def\math{\fam\mathfam}
     \def \\ { \cr }
     \def\R{{\math R}}
     
     \def\N{{\math N}}
     \def\E{{\math E}}
     \def\P{{\math P}}

     \def\D{{\math D}}

     \def \e{{\rm e}}

     \def \f{{\cal F}}

     \def \k{{\tt k}}

\def\la{\longrightarrow}
\def\ba{\begin{eqnarray*}}
\def\ea{\end{eqnarray*}}

\def\build#1_#2^#3{\mathrel{\mathop{\kern 0pt#1}\limits_{#2}^{#3}}}
\def\be{\begin{equation}}
\def\ee{\end{equation}}
     \newtheorem{theorem}{Theorem}
     \newtheorem{proposition}{Proposition}
     \newtheorem{lemma}{Lemma}

\begin{document}

     \centerline{\LARGE \bf Stochastic flows associated to}
     \vskip 2mm
     \centerline{\LARGE \bf  coalescent processes III:}
     \vskip 2mm
     \centerline{\LARGE \bf  Limit theorems}

     \vskip 1cm

\centerline{\bf Dedicated to the memory of Joseph Leo DOOB}

\vskip 1cm
     \centerline{\Large \bf Jean Bertoin$^{(1)}$ and Jean-Fran\c cois Le
     Gall$^{(2)}$}
     \vskip 1cm
     \noindent
     \noindent
     (1) {\sl Laboratoire de Probabilit\'es et Mod\` eles Al\'eatoires
     and Institut universitaire de France,
        Universit\'e Pierre et Marie Curie,  175, rue du Chevaleret,
        F-75013 Paris, France.}
     \vskip 2mm
     \noindent
     (2) {\sl DMA, Ecole normale sup\'erieure, 45, rue d'Ulm, F-75005 Paris,
     France.}
     \vskip 15mm

     \noindent{\bf Summary. }{\small
We prove several limit theorems that relate coalescent processes to
continuous-state
branching processes. Some of these theorems are stated in terms of the
so-called generalized Fleming-Viot processes, which describe
the evolution of a population with fixed size, and are duals to the
coalescents with multiple collisions studied
by Pitman and others. We first discuss asymptotics when the initial
size of the population tends to infinity.
In that setting,  under appropriate hypotheses,
  we show that a rescaled version of the generalized Fleming-Viot
process converges weakly to
a continuous-state branching process. As a corollary, we get a
hydrodynamic limit for certain sequences of
coalescents with multiple collisions: Under an appropriate scaling, the
empirical measure associated with sizes of the blocks  converges to
a (deterministic) limit which solves a generalized form of
Smoluchowski's coagulation equation.
We also study the behavior in small time of a fixed coalescent with
multiple collisions, under
a regular variation assumption on the tail of the measure $\nu$
governing the coalescence events. Precisely, we prove that
the number of blocks with size less than $\varepsilon x$ at time
$(\varepsilon\nu([\varepsilon,1]))^{-1}$ behaves like
$\varepsilon^{-1}\lambda_1(]0,x[)$ as $\varepsilon\to 0$, where
$\lambda_1$ is the distribution of the size of one cluster at time $1$
in a continuous-state branching process with stable branching mechanism. This
generalizes a classical result for the Kingman coalescent.}

   \vskip 3mm
     \noindent
        {\bf Key words. \small Flow, coalescence, Fleming-Viot process,
continuous state branching process, Smoluchowski's coagulation equation.}
        \vskip 5mm
     \noindent
     {\bf A.M.S. Classification.}  {\tt 60 G 09, 60 J 25, 92 D 30.}
     \vskip 3mm
     \noindent{\bf e-mail.} {\tt $(1):$ jbe@ccr.jussieu.fr , $(2):$
     legall@dma.ens.fr}

     \section{Introduction}
     J.~L.~Doob was a
pioneer in the development of the theory of martingales and its applications to
probability theory, potential theory or functional analysis.
The fundamental contributions that he made in this field
form the cornerstones of one of the richest
veins explored in mathematics during the last half-century.
In particular, martingales and stochastic calculus provide nowadays
key tools for
studying the asymptotic behavior of random processes; see the classical books
by Ethier and Kurtz \cite{EK} and Jacod and Shiryaev \cite{JS}.
In the present work, we shall apply such techniques to investigate a class of
stochastic flows related to certain population dynamics.

The general motivation for the present work is to get a better understanding
of the relations between the so-called coalescents with multiple collisions,
which were introduced independently by Pitman \cite{Pi}
and Sagitov \cite{Sa}, and continuous-state branching processes.
Note from \cite{Sa} that coalescents with multiple collisions can be viewed as
asymptotic models for the genealogy of a discrete population with a
fixed size, and so
the existence of connections with branching processes should not come
as a surprise.
Such connections were already derived in \cite{BeLG0}, where the
Bolthausen-Sznitman
coalescent was shown to describe the genealogical structure of a particular
continuous-state branching process introduced by Neveu, and in
\cite{7mercenaires},
where similar relations were obtained between the so-called beta-coalescents
and continuous-state branching processes with stable branching mechanism.
Here, we do not focus on exact distributional identities, but
rather on asymptotics for functionals of coalescent processes, where the
limiting objects are given in terms of branching processes. In order
to get such
asymptotics, we apply the machinery of limit theorems
for semimartingales \cite{JS} to the so-called generalized Fleming-Viot
processes, which where shown in \cite{BeLG1} to be duals to the coalescents
with multiple collisions.

Generalized Fleming-Viot processes, which model the evolution of a
continuous population with fixed size $1$,
have appeared in articles  by Donnelly and Kurtz \cite{DK1, DK2}, and
were studied more recently
in our work \cite{BeLG1, BeLG2}.  It is convenient
to view a generalized Fleming-Viot process as a stochastic flow
$(F_t, t\geq0)$ on $[0,1]$, such that for each
$t\geq0$, $F_t: [0,1]\to[0,1]$ is a (random) right-continuous
increasing map with
$F_t(0)=0$ and $F_t(1)=1$.
We should think of  the unit interval
as a population, and then of $F_{t}$ as the distribution function
of a (random) probability measure $dF_{t}(x)$ on $[0,1]$. The
evolution of the latter
is related to the dynamics of the population as follows :
For every $0\leq r_1< r_2 \leq 1$, the
interval $]F_t(r_1), F_t(r_2)]$
represents  the sub-population at time $t$ which
consists of descendants of the sub-population $]r_1,r_2]$ at the initial time.
The transitions of the flow are Markovian, and more precisely, for
every $s,t\geq 0$,
we have $F_{t+s}=\tilde F_s\circ F_t$, where $\tilde F_s$ is a copy
of $F_s$ independent of $(F_r,0\leq r\leq t)$. The distribution of the flow is
then characterized  by a measure $\nu$ on $]0,1]$ such that
$\int_{]0,1]}x^2\nu(dx)<\infty$. To explain this, consider the simple
case where
$\nu$ is a finite measure. Let $((T_i,U_i,\xi_i), i\in\N)$
denote the sequence of
atoms of a Poisson random measure on $[0,\infty[\times [0,1]\times [0,1]$
with intensity $dt\otimes du\otimes \nu(dx)$, ranked in the
increasing order of the first coordinate.
The process $(F_t, t\geq0)$ starts from
$F_0={\rm Id}$, remains constant on the intervals
$[T_{i-1}, T_i[$ (with the usual convention that $T_0=0$), and for every
$i\in\N$
$$F_{T_i}\,=\, \Delta_i\circ F_{T_{i-1}}$$
where
$$\Delta_i(r)=\xi_i{\bf 1}_{\{U_i\leq r\}}+r(1-\xi_i)\,,\qquad
r\in[0,1].$$ In terms of the population model, this means that at
each time $T_i$,
an individual
in the population at time $T_{i-1}$ is picked uniformly at random and
gives birth to a sub-population of size $\xi_i$. Simultaneously, the rest
of the population shrinks by factor $1-\xi_i$, so the total size of the
population remains $1$. The previous description does not apply when
$\nu$ is infinite, since then the Poisson measure will have infinitely
many atoms on a finite time interval. Still, the Fleming-Viot flow
can be constructed
via a suitable limiting procedure (\cite{BeLG1} Theorem 2).

Our first motivation for studying generalized Fleming-Viot processes
came from their
remarkable connection \cite{BeLG1} with the class of
coalescents with multiple collisions considered by Pitman \cite{Pi} and
Sagitov \cite{Sa}. To describe this connection, fix some time $T>0$
viewed as the present date at which the population is observed,
and pick a sequence of individuals labelled $1,2,\ldots$ independently and
uniformly over $[0,1]$.
For every $t\leq T$,
we obtain a partition $\Pi(t)$ of $\N$ by
gathering individuals having the
same ancestor at time $T-t$. The process $(\Pi(t),0 \leq t \leq T)$
is then a Markovian coalescent process on the space of partitions of $\N$.
In the terminology of \cite{Pi},
it is a $\Lambda$-coalescent, with $\Lambda(dx)=x^{2}\nu(dx)$,
started  from the partition of
$\N$ into singletons.
As a consequence of Kingman's theory of exchangeable partitions,
for every $t\geq 0$, each
block of $\Pi(t)$ has an asymptotic frequency, also called the size
of the block,
and the ranked sequence of these
frequencies yields a Markov process called the mass-coalescent. As a
consequence of
the preceding construction, the mass-coalescent at time $t$ has the same
distribution as the ranked sequence of jump sizes of $F_t$.

The first purpose of the present work is to investigate the
asymptotic behavior of a rescaled version of the preceding population model.
Specifically,  we consider a family $(\tilde\nu^{(a)}, a>0)$ of measures
on $]0,1]$ such that $\int_{]0,1]}x^2\tilde\nu^{(a)}(dx)<\infty$ for
every $a>0$, and the associated
generalized Fleming-Viot processes $\tilde F^{(a)}$. For each $a>0$, we rescale
$\tilde F^{(a)}$ by a factor $a$ in space and time, i.e. we set
$$F^{(a)}_t(r)\,:=\,a\tilde F^{(a)}_{at}(r/a)\,,\qquad r\in[0,a],
t\geq0\,.$$
So the process $F^{(a)}$ describes the evolution of a population
with fixed size $a$. Roughly speaking, considering $F^{(a)}$
in place of $\tilde F^{(a)}$ enables us to focus on the dynamics
of a sub-population having size of order $1/a$.
Denote by $\nu^{(a)}$ the image of
$\tilde\nu^{(a)}$ under the dilation $x\to ax$, and assume that the measures
$(x^2\wedge x)\nu^{(a)}(dx)$ converge weakly as $a\to\infty$ to a
finite measure on
$]0,\infty[$, which we may write in the form $(x^2\wedge x)\pi(dx)$.
Then Theorem \ref{largepop} shows that $F^{(a)}$ converges in distribution
to the critical continuous-state branching process $Z$ with branching
mechanism
$$\Psi(q)\,=\,\int_{]0,\infty[}(\e^{-qx}-1+qx)\pi(dx)\,,\qquad q
\geq0\,.$$

As a consequence of this limit theorem, we derive a hydrodynamic
limit for the associated coalescent
processes (Theorem \ref{T3}).
Precisely,  we show that under the same assumptions as above, for
every $t\geq 0$, the
empirical measure corresponding to the jumps of $\tilde F^{(a)}_t$
(or equivalently to
the block sizes in the associated coalescent) converges, modulo a
suitable rescaling,
towards a deterministic measure $\lambda_t$. Informally, $\lambda_t$
is the distribution ofa cluster at time
$t$, that is a collection of individuals sharing the same ancestor at
the initial time, in the
continuous-state branching process with branching mechanism $\Psi$.
In a way analogous
to the derivation of Smoluchovski's coagulation equation from
stochastic models (see Aldous \cite{Al}, Norris
\cite{Norris} and the references therein for background) we prove
that the family $(\lambda_t,t>0)$ solves
a generalized coagulation equation of the form
$${d \langle\lambda_t,f\rangle\over dt}=
\sum_{k=2}^\infty {(-1)^k\Psi^{(k)}(\langle \lambda_t,1\rangle)\over k!}
\int_{]0,\infty[^k}
(f(x_1+\cdots+x_k)-(f(x_1)+\cdots+f(x_k)))\,\lambda_t(dx_1)\ldots
\lambda_t(dx_k)$$
where $f$ can be any continuous function with compact support on
$]0,\infty[$ (Proposition \ref{Smolu}).

In the last part of this work, we study the small time behavior of generalized
Fleming-Viot processes and $\Lambda$-coalescents, under a regular variation
assumption on the measure $\nu$ (recall that
$\Lambda(dx)=x^2\nu(dx)$). Precisely, we
assume that the tail $\nu([\varepsilon,1])$ is regularly varying with
index $-\gamma$
when $\varepsilon$ goes to $0$.
We are interested in the case when the $\Lambda$-coalescent comes
down from infinity
(i.e. for every $t>0$, $\Pi_t$ has finitely many blocks), which
forces $1\leq \gamma\leq 2$.
Leaving aside the boundary cases we suppose that $1<\gamma<2$. As a
consequence of Theorem \ref{largepop},
we prove that the rescaled Fleming-Viot process
$$F^\varepsilon_t(x):={1\over
\varepsilon}\,F_{t/(\varepsilon\nu([\varepsilon,1]))}(\varepsilon x)$$
converges in distribution to the continuous-state branching process with stable
branching mechanism:
$$\Psi_\gamma(q)={\Gamma(2-\gamma)\over \gamma-1}\,q^\gamma.$$
We then use this result to investigate the small time behavior of the
size of blocks in the $\Lambda$-coalescent. Write $N_t(]0,x[)$ for
the number of blocks with size
less than $x$ in the $\Lambda$-coalescent at time $t$. If
$g(\varepsilon)=(\varepsilon
\nu([\varepsilon,1]))^{-1}$, Theorem \ref{blocks} states that
$$\sup_{x\in]0,\infty[}\Big|\varepsilon
N_{g(\varepsilon)}(]0,\varepsilon x[) -\lambda_1 (]0,x[)\Big|
\build{\la}_{\varepsilon\to 0}^{} 0,$$
in probability.
Furthermore, the measure $\lambda_1$ can be characterized by its
Laplace transform
$$\int
(1-\e^{-qr})\lambda_1(dr)=(\Gamma(2-\gamma)+q^{1-\gamma})^{1/(1-\gamma)}.$$
Theorem \ref{blocks} is analogous to a classical result for the sizes
of blocks in
the Kingman coalescent in small time (see Aldous \cite{Al}). The proof uses
an intermediate estimate for the total number of blocks in a
$\Lambda$-coalescent, which is closely related to the recent paper
\cite{BBS} dealing with beta-coalescents.

The paper is organized as follows. Section 2 gives a few preliminary
results about
continuous-state branching processes. In particular, the Poisson
representation (Proposition \ref{T1})
may have other applications. Section 3 states our first limit theorem
for generalized
Fleming-Viot processes. The derivation of the hydrodynamic limit is
developed in Section 4,
which also discusses the generalized coagulation equation for the
family $(\lambda_t,t\geq 0)$.
Finally Section 5 is devoted to the behavior in small time of
generalized Fleming-Viot processes
and $\Lambda$-coalescents.

\noindent{\bf Notation}. We use the notation $\langle\mu,f\rangle$
for the integral of the function $f$
with respect to the measure $\mu$. We denote by ${\cal M}_{\rm F}$ the
space of all finite measures on
$]0,\infty[$, which is equipped with the usual weak topology. We also
denote by ${\cal M}_{\rm R}$ the space of all Radon
measures on $]0,\infty[$. The set ${\cal M}_{\rm R}$ is equipped with
the vague topology: A sequence $(\mu_n,n\in\N)$
in  ${\cal M}_{\rm R}$ converges to $\mu\in {\cal M}_{\rm R}$
if and only if for every continuous function $f:]0,\infty[\to\R$ with compact
support, $\lim_{n\to\infty}\langle \mu_n,f\rangle=\langle \mu,f\rangle$.

\section{Stochastic flows of branching processes}

In this section, we give a few properties of
continuous-state branching processes that will
be needed in the proof of our limit theorems.
A {\it critical branching mechanism} is a function
$\Psi:[0,\infty[\rightarrow [0,\infty[$ of the type
\be
\label{mechbranch}
\Psi(q)\,=\,\beta q^2 + \int_{]0,\infty[}\left(\e^{-rq}
-1+rq\right)\pi(dr)
\ee
where  $\beta\geq0$ is the so-called Gaussian coefficient and $\pi$ is a
measure on
$]0,\infty[$ such that $\int
(r\wedge r^2)\pi(dr)<\infty$. The continuous-state branching process
with branching mechanism $\Psi$ (in short the $\Psi$-CSBP) is the Markov
process with values in $\R_+$, whose transition kernels $Q_t(x,dy)$ are
determined by the  Laplace transform
\be
\label{Laplace}
\int Q_t(x,dy)\,\e^{-q y}=\exp(-x\,u_t(q))\;,\ x,t\geq
0,\;q\geq 0\,,
\ee
where the function
$u_t(q)$ solves
\begin{equation}\label{eqsg}
{\partial u_t(q)\over \partial t}\,=\,-\Psi(u_t(q))\quad ,\quad
u_0(q)\,=\,q\,.\end{equation}
The criticality of $\Psi$ implies that a $\Psi$-CSBP is a
nonnegative martingale.
If $Z^1$ and $Z^2$ are two independent $\Psi$-CSBP's started respectively
at $x_1$ and $x_2$, then $Z^1+Z^2$ is also a $\Psi$-CSBP, obviously with
initial value $x_1+x_2$. From this additivity or branching property,
we may construct a two-parameter process $Z=(Z(t,x),t,x\geq 0)$, such
that:
\begin{description}
\item{$\bullet$} For each fixed $x\geq 0$, $(Z(t,x),t\geq 0)$ is a
$\Psi$-CSBP
with c\` adl\` ag paths and initial value $Z(0,x)=x$.
\item{$\bullet$} If $x_1,x_2\geq 0$, $Z(\cdot, x_1+x_2)-Z(\cdot,x_1)$ is
independent of the processes
$\left(Z(\cdot,x), 0\leq x
\leq x_1\right)$ and has the same law as $Z(\cdot,x_2)$.
\end{description}

These properties entail that for each fixed $t\geq0$, $Z(t,\cdot)$ is an
increasing process with independent and stationary increments.
Its right-continuous version is a subordinator with Laplace exponent $u_t$
determined by
(\ref{Laplace}) and (\ref{eqsg}).
By the L\'evy-Khintchin formula, there exists a unique drift
coefficient $d_t\geq0$ and a unique measure $\lambda_t$ on $]0,\infty[$
with $\int_{]0,\infty[}(1\wedge x)\lambda_t(dx)<\infty$ such that
\be
\label{driftLevy}
u_t(q)\,=\,qd_t+\int_{]0,\infty[}(1-\e^{-qx})\lambda_t(dx)\,,\qquad
q\geq0\,.
\ee
One  refers to $\lambda_t$ as the L\'evy measure of $Z(t,\cdot)$.
Measures $\lambda_t$ play an important role in this work. Informally,
we may say that $\lambda_t$ is the `distribution' of the size of
the set of descendants at time $t$ of a single individual at time $0$.
This assertion is informal since $\lambda_t$ is not a probability
distribution (it may even be an infinite measure). A correct way of
stating the above (in the case $d_t=0$) is as follows: $Z(t,x)$
is the sum of the atoms of a Poisson measure with intensity
$x\lambda_t(\cdot)$.
Moreover, the study of the genealogical structure of the $\Psi$-CSBP (see e.g.
\cite{DuLG}) allows one
to interpret each of these atoms as the size of a family of
individuals at time $t$ that have the same ancestor at the initial time.

 From now on, we assume that $\beta=0$ and we exclude the trivial case
$\pi=0$.

We start by recalling in our special case an important connection between
continuous-state branching processes and L\'evy processes due to
Lamperti \cite{La}.
Let $x>0$ be fixed, and let $\xi=(\xi_t, t\geq0)$ denote a real-valued
L\'evy process with no negative jumps, started from $\xi_0=x$, and whose
Laplace exponent is specified by
$$\E\left[\exp (-q(\xi_t-\xi_0))\right]\,=\,\exp t\Psi(q)\,,\qquad
q\geq0\,.$$
In particular $\pi$ is the L\'evy measure of $\xi$.
The criticality of the branching mechanism $\Psi$ ensures that
the L\'evy process $\xi$ has centered increments and thus oscillates.
In particular the first passage time
  $\zeta:=\inf\left\{t\geq0: \xi_t=0\right\}$ is finite a.s. Next,
introduce for every $t\geq0$
$$\gamma(t)\,=\,\int_{0}^{t\wedge \zeta}{ds\over \xi_s}\quad ,\quad
C_t\,=\,\inf\left\{s\geq0: \gamma(s)>t\right\}\wedge \zeta\,.$$
Then the time-changed process
$\left(\xi\circ C_t, t\geq0\right)$ has the same distribution as
$\left(Z(t,x), t\geq0\right)$.

It follows from this representation that $(Z(t,x),t\geq 0)$ is
a purely discontinuous martingale. We can also use the Lamperti
transformation to calculate the compensator of the jump measure
of this martingale.
By the L\'evy-It\^o decomposition, the compensator of the jump measure
of $\xi$,
$$\sum_{\{t:\Delta\xi_t\neq0\}}\delta_{(t,\Delta\xi_t)},$$
is $dt\otimes \pi(dx)$. By a time change argument, we can then deduce that
the compensator of the measure
$$\sum_{\{t:\Delta Z(t,x)\neq 0\}}\delta_{(t,\Delta Z(t,x))}$$
is $Z(t,x)dt\otimes \pi(dr)$.

Since $(Z(t,x),t\geq 0)$ is a purely discontinuous martingale, the
knowledge of the compensator of its jump measure completely
determines the characteristics of this semimartingale, in the sense
of \cite{JS} Chapter II. We will need the fact that the distribution
of $(Z(t,x),t\geq 0)$, and more generally of the
multidimensional process $((Z(t,x_1),Z(t,x_2),\ldots,Z(t,x_p));t\geq 0)$
for any choice of $p$ and $x_1,\ldots,x_p$, is
uniquely determined by its characteristics.

Fix an integer $p\geq1$ and
define
\begin{equation}
\label{eqdomain}
{\cal D}_p:=\{x=(x_1,\ldots,x_p):0\leq x_1\leq x_2\leq \cdots\leq
x_p\}.
\end{equation}
For every $(y_1,\ldots,y_p)\in{\cal D}_p$, define a $\sigma$-finite measure
$U(y_1,\ldots,y_p;dz_1,\ldots,dz_p)$ on $\R_+^p\backslash\{0\}$
by setting, for any measurable function $\varphi:\R_+^p\to\R_+$ that
vanishes at
$0$,
\be
\label{jumpCSBP}
\int
U(y_1,\ldots,y_p;dz_1,\ldots,dz_p)\,\varphi(z_1,\ldots,z_p)
=\int\pi(dr)\int_0^\infty du\,\varphi(r{\bf 1}_{\{u\leq y_1\}},\ldots,r{\bf
1}_{\{u\leq y_p\}}).
\ee

\begin{proposition}
\label{UniCSBP}
Let $(x_1,\ldots,x_p)\in{\cal D}_p$ and let
$(Z^1,\ldots,Z^p)$ be a $p$-dimensional semimartingale taking values
in ${\cal D}_p$, such that $(Z^1_0,\ldots,Z^p_0)=(x_1,\ldots,x_p)$.
The following two properties are equivalent:
\begin{description}
\item{\rm (i)} The processes $((Z^1_t,\ldots,Z^p_t);t\geq 0)$
and $((Z(t,x_1),\ldots,Z(t,x_p));t\geq 0)$
have the same distribution.
\item{\rm(ii)} The process $((Z^1_t,\ldots,Z^p_t);t\geq 0)$ is a purely
discontinuous local martingale, and the compensator of its jump measure is
the measure
$$\theta(dt,dz_1\ldots dz_p)=dt\;U(Z^1_t,\ldots,Z^p_t;dz_1,\ldots,dz_p).$$
\end{description}
\end{proposition}

\proof The implication (i)$\Rightarrow$(ii) is a straightforward
consequence  of the remarks preceding the statement and the branching
property of continuous-state branching processes. We concentrate on the
proof of the converse implication (ii)$\Rightarrow$(i). Let
$q=(q_1,\ldots,q_p)\in]0,\infty[^p$, and let $Y_t=(Y^1_t,\ldots,Y^p_t)$
be defined by $Y^i_t=Z^i_t-Z^{i-1}_t$ if $i\geq 2$ and $Y^1_t=Z^1_t$.
Notice that $Y^i_t\geq 0$. Using property (ii), an application
of It\^o's formula (cf Theorem II.2.42 in \cite{JS}) yields that
the process
\begin{eqnarray*}
& &\exp(- q\cdot Y_t)-\exp(- q\cdot Y_0)\\
&&-\sum_{i=1}^{p}\int_{[0,t]\times [0,\infty[\times ]0,\infty[}
\exp(-q\cdot Y_s)\left(\e^{-q_i r}-1+q_i r\right){\bf 1}_{
\{u\leq Y^i_s\}}ds\,du\,\pi(dr)
\end{eqnarray*}
is a local martingale. This local martingale is bounded
over the time interval $[0,t]$ for any $t\geq 0$, hence is a martingale.
Taking expectations leads to
\be
\label{unitech}
\E[\e^{-q\cdot Y_t}]=\E[\e^{-q\cdot Y_0}]
+\sum_{i=1}^p \Psi(q_i)\int_0^t ds\,\E[Y^i_s\,\e^{-q\cdot Y_s}].
\ee
It is immediate to verify from (ii) that each $Y^i$ is also
a nonnegative local martingale, and so $\E[Y^i_s]\leq \E[Y^i_0]=x_i-x_{i-1}$
(by convention $x_0=0$). If we set
$f_t(q)=\E[\e^{-q\cdot Y_t}]$
we have
$${\partial f_t(q)\over \partial q_i}=-\E[Y^i_t\,\e^{-q\cdot Y_t}]$$
and so we deduce from (\ref{unitech}) that
\be
\label{unitech2}
{\partial f_t(q)\over \partial t}
+
\Psi(q)\cdot \nabla f_t(q)\,=\,0\,,
\ee
where we write $\Psi(q)=(\Psi(q_1),\ldots,\Psi(q_p))$. In order
to solve (\ref{unitech2}), fix $t_1>0$, and consider the function
$g(t)=(u_{t_1-t}(q_1),\ldots,u_{t_1-t}(q_p))$ for $t\in[0,t_1]$,
where $u_t(q)$ is as in (\ref{eqsg}).
Since $$g'(t)=(\Psi(u_{t_1-t}(q_1)),\ldots,\Psi(u_{t_1-t}(q_p))),$$
it follows that
$${\partial f_t\circ g\over \partial t}={\partial f_t\over \partial t}
\circ g+ g'(t)\cdot \nabla f_t(g(t))=0$$
by (\ref{unitech2}). Hence $f_t\circ g(t)$ is constant over $[0,t_1]$, and
$$f_{t_1}(q)=f_{t_1}(g(t_1))=f_0(g(0))=\exp(-\sum_{i=1}^p (x_i-x_{i-1})
u_{t_1}(q_i)).$$
This shows that
$$(Y^1_{t_1},\ldots,Y^p_{t_1})
\build=_{}^{\rm(d)}(Z(t_1,x_1),Z(t_1,x_2)-Z(t_1,x_1),
\ldots,Z(t_1,x_p)-Z(t_1,x_{p-1}))$$
and so
$$(Z^1_{t_1},\ldots,Z^p_{t_1})\build=_{}^{\rm(d)}
(Z(t_1,x_1),Z(t_1,x_2),\ldots,Z(t_1,x_p)).$$
It is easy to iterate this argument to obtain that the processes
$((Z^1_t,\ldots,Z^p_t);t\geq 0)$ and $((Z(t,x_1),\ldots,Z(t,x_p));t\geq 0)$
have the same finite-dimensional marginal distributions. The desired result
follows since both processes have c\` adl\` ag paths. \QED

We now turn our attention to the representation of critical CSBP as
stochastic flows on
$[0,\infty[$ solving simple stochastic differential equations.
On a suitable
filtered probability space $(\Omega,\f,(\f_t),\P)$, we consider :

\noindent $\bullet$  an $(\f_t)$-Poisson random measure
$$M\,=\,\sum_{i=1}^{\infty}\delta_{(t_i, u_i, r_i)}\,,$$
on
$\R_+\times [0,\infty[\times
]0,\infty[$, with intensity $dt\otimes du\otimes \pi(dr)$.

\noindent $\bullet$ a collection $(X_t(x), t\geq0)$, $x\in\R_+$ of
c\` adl\` ag $(\f_t)$-martingales with values in $\R_+$,

\noindent $\bullet$ the stochastic differential equation
\begin{equation}\label{eq2}
X_t(x)\,=\,x+ \int_{[0,t]\times
[0,\infty[\times]0,\infty[}M(ds,du,dr)\,r\,{\bf 1}_{\{u\leq
X_{s-}(x)\}}\,.
\end{equation}
The Poissonian stochastic integral in the right-hand side should be
understood with respect to the compensated Poisson measure $M$ (see
e.g. Section II.1 of \cite{JS}). This stochastic integral
is well defined according to Definition II.1.37 of \cite{JS},
since the increasing process
$$t\longrightarrow \Big(\int_{[0,t]\times [0,\infty[\times
[0,\infty[} M(ds,du,dr)\,r^2\,{\bf 1}_{\{u\leq
X_{s-}(x)\}}\Big)^{1/2}$$
is locally integrable under our assumption on $\pi$.

A pair $\left(M,(X_{\cdot}(a), a\geq0)\right)$ satisfying the
above conditions
will be called a {\it weak solution} of (\ref{eq2}).

\begin{proposition}\label{T1}  The equation {\rm(\ref{eq2})} has a weak
solution
  which satisfies the
additional property that $X_t(x_1)\leq X_t(x_2)$ for every $t\geq 0$, a.s.
whenever $0\leq x_1\leq x_2$. Moreover, for
every such solution $(M,X)$,  for every $p\in \N$ and $0\leq x_1\leq
\ldots\leq x_p$, the process
$((X_t(x_1),\ldots,X_t(x_p)), t\geq0)$
has the same distribution as $((Z(t,x_1),\ldots,Z(t,x_p)), t\geq0)$.
\end{proposition}

\proof The second part of the statement is immediate from the
implication (ii)$\Rightarrow$(i) in Proposition \ref{UniCSBP}.
The first part can be deduced from Theorem 14.80 in \cite{Jac}
by the same arguments that were used in the proof of
Theorem 2 in \cite{BeLG2}. We leave details to the reader as this result is
not really needed below except for motivation.
\QED

\section{Generalized Fleming-Viot flows and their limits}
We now recall some results from \cite{BeLG1, BeLG2} on generalized
Fleming-Viot processes and related stochastic flows. Let  $\nu$ denote a
$\sigma$-finite measure on
$]0,1]$ such that $\int_{]0,1]}x^2\nu(dx)<\infty$. According to Section
5.1 in \cite{BeLG1}, one can associate with $\nu$ a Feller
process
$(F_t, t\geq0)$ with values in the space of distribution functions of
probability measures on
$]0,1]$ (i.e. for each $t\geq0$, $F_t$ is a c\` adlag increasing map
from $[0,1]$ to $[0,1]$ with $F_t(0)=0$ and $F_t(1)=1$),
whose evolution is characterized by $\nu$ and has been described in Section 1
in the special case when $\nu$ is finite.

In \cite{BeLG2}, we have shown that such generalized Fleming-Viot
processes can be described as the solution to a certain system of
Poissonian SDE's. More precisely, on a suitable
filtered probability space $(\Omega,\f,(\f_t),\P)$, one can construct the
following processes:

\noindent $\bullet$  an $(\f_t)$-Poisson point process $N$ on
$\R_+\times ]0,1[\times
]0,1]$ with intensity $dt\otimes du \otimes
\nu(dr)$,

\noindent $\bullet$ a collection $(Y_t(x), t\geq0)$, $x\in[0,1]$, of
adapted c\` adl\` ag processes with values in $[0,1]$
with $Y_t(x_1)\leq Y_t(x_2)$ for all $t\geq0$ a.s. when $0\leq x_1 \leq
x_2 \leq 1$,

\noindent in such a way that for every $r\in[0,1]$, a.s.
\begin{equation}\label{EQ21}
Y_t(x)\,=\,x+ \int_{[0,t]\times
]0,1[\times]0,1]}N(ds,du,dr)\,r\,\left({\bf 1}_{\{u\leq Y_{s-}(x)\}}
-Y_{s-}(x)\right)\,.
\end{equation}
The Poissonian stochastic integral in the right-hand side should
again be
understood with respect to the compensated Poisson measure $N$.

Weak uniqueness holds for this
system of SDE's (Theorem 2 in \cite{BeLG2}). Furthermore,
for every integer $p\geq 1$ and every $0\leq x_1\leq \ldots \leq x_p\leq
1$,
the processes
$((Y_t(x_1),\ldots,Y_t(x_p)), t\geq0)$ and $((F_t(x_1),\ldots,F_t(x_p)),
t\geq0)$ have the same distribution. Note the similarity with Proposition
\ref{T1}: Compare (\ref{eq2}) and (\ref{EQ21}). This strongly suggests to
look for asymptotic results relating the processes $Z(t,x)$ and $F_t(x)$.

For every integer $p\geq 1$ and every $a>0$, set
$${\cal
D}^{a}_p={\cal D}_p\cap [0,a]^p
=\{(x_1,\ldots,x_p)\in\R_+^p:0\leq x_1\leq x_2\leq \ldots\leq x_p\leq a\}.$$
 From (\ref{EQ21}) we see that for every $(x_1,\ldots,x_p)\in{\cal
D}^1_p$, the process $(F_t(x_1),\ldots,F_t(x_p))$ is a purely discontinuous
martingale, and the compensator of its jump measure is
$$dt\;R(F_t(x_1),\ldots,F_t(x_p);dz_1,\ldots,dz_p),$$
where  for every $(y_1,\ldots,y_p)\in{\cal D}^1_p$, the measure
$R(y_1,\ldots,y_p;dz_1,\ldots,dz_p)$ on $\R^p\backslash\{0\}$ is
determined by
$$\int R(y_1,\ldots,y_p;dz_1,\ldots,dz_p)\,\varphi(z_1,\ldots,z_p)
=\int\nu(dr)\int_0^1 du\,\varphi(r({\bf
1}_{\{u\leq y_1\}}-y_1),\ldots,r({\bf
1}_{\{u\leq y_p\}}-y_p)).$$

Consider now a family $(\tilde \nu^{(a)}, a>0)$ of measures on $]0,1]$
with $\int_{]0,1]}r^2\tilde \nu^{(a)}(dr)<\infty$, and for each $a>0$,
let $\tilde F^{(a)}$ be the associated Fleming-Viot process.
We then write
$$F^{(a)}_t(x):=a\tilde  F^{(a)}_{at}(x/a)\,,\qquad x\in[0,a], t\geq0$$
for the rescaled version of the Fleming-Viot flow.
So, for each $t\geq0$, $F^{(a)}_t$ is the distribution function
of a measure on $]0,a]$ with total mass $a$.
For every fixed
real number $a>0$, we also denote by $\nu^{(a)}$ the measure on $]0,\infty[$
which is $0$ on $]a,\infty[$ and whose restriction to $]0,a]$ is given by
the image of $\tilde \nu^{(a)}$
under the dilation
$r\to ar$ from $]0,1]$ to $]0,a]$. In particular
$r^2\nu^{(a)}(dr)$ is a finite measure on $]0,\infty[$.

By a scaling argument, we see that,
for every $(x_1,\ldots,x_p)\in{\cal D}^a_p$, $(F^{(a)}_t(x_1),\ldots,
F^{(a)}_t(x_p))$ is a purely discontinuous martingale,
with values in ${\cal D}^a_p$,
and the compensator of its jump measure
is
\be
\label{compensFV}
\mu_{(a)}(dt,dz_1\ldots
dz_p)=dt\,R^{(a)}(F^{(a)}_t(x_1),\ldots,F^{(a)}_t(x_p); dz_1,\ldots,dz_p)
\ee
where
\begin{eqnarray}
\label{jumpFVbis}
&&\int
R^{(a)}(y_1,\ldots,y_p;dz_1,\ldots,dz_p)\,\varphi(z_1,\ldots,z_p)
\nonumber\\
&&\quad =\int\nu^{(a)}(dr)\int_0^{a} du\,\varphi(r({\bf
1}_{\{u\leq y_1\}}-a^{-1} y_1),\ldots,r({\bf 1}_{\{u\leq
y_p\}}-a^{-1} y_p)).
\end{eqnarray}

Let $\pi$ be as in Section 2 a nontrivial measure on $]0,\infty[$ such that
$\int (r\wedge r^2)\pi(dr)<\infty$, and let $\Psi$ be as in
(\ref{mechbranch}). Denote by $(Z(t,x),t\geq 0,x\geq 0)$
the associated flow of continuous-state branching processes constructed
in Section 2.

\noindent{\bf Assumption} (H) {\it The measures $(r\wedge r^2)
\nu^{(a)}(dr)$ converge to $(r\wedge r^2)\pi(dr)$ as $a\to\infty$,
in the sense of weak convergence in ${\cal M}_{\rm F}$.}

\begin{theorem}
\label{largepop}
Under Assumption {\rm (H)},
for every $(x_1,\ldots,x_p)\in{\cal D}_p$,
$$((F^{(a)}_t(x_1),\ldots,F^{(a)}_t(x_p));t\geq 0)
\build{\la}_{a\to \infty}^{\rm(d)} ((Z(t,x_1),\ldots,Z(t,x_p));
t\geq 0)$$
in the Skorokhod space $\D(\R_+,\R^p)$.
\end{theorem}

\proof The proof only uses the facts that $(F^{(a)}_t(x_1),\ldots,
F^{(a)}_t(x_p))$ is a purely discontinuous martingale and that
the compensator of its jump measure is given by (\ref{compensFV})
and (\ref{jumpFVbis}). The latter properties indeed characterize
the law of the process $(F^{(a)}_t(x_1),\ldots,
F^{(a)}_t(x_p))$ (cf Lemma 1 in \cite{BeLG2}), but we do not
use this uniqueness property in the proof. We fix a sequence
$(a_n)$ tending to $+\infty$, and
$(x_1,\ldots,x_p)\in{\cal D}_p$. To simplify notation we write
$$Y^n_t=(Y^{n,1}_t,\ldots,Y^{n,p}_t)=
(F^{(a_n)}_t(x_1),\ldots,F^{(a_n)}_t(x_p))$$
which makes sense as soon as $a_n\geq x_p$, hence for all $n$
sufficiently large. We also set
$$Z_t=(Z^1_t,\ldots,Z^p_t)=(Z(t,x_1),\ldots,Z(t,x_p)).$$

We rely on general limit theorems for semimartingales with jumps
which can be found in the book \cite{JS}. To this end, we first need
to introduce a truncation function $h:\R\to \R$, that is a
bounded continuous function such that $h(x)=x$ for every
$x\in[-\delta,\delta]$, for some $\delta>0$. We may and will assume
that $h$ is nondecreasing, $|h(x)|\leq |x|\wedge 1$ for every $x\in \R$
and that $h$ is Lipschitz continuous with Lipschitz constant $1$,
that is $|h(x)-h(y)|\leq |x-y|$ for every $x,y\in\R$. We can then
consider the associated (modified) triplet of characteristics
of the $p$-dimensional semimartingale
$Y^n$:
$$(B^n,\tilde C^n,\mu_{(a_n)}).$$
See Definition II.2.16 in \cite{JS}. To be specific,
$\mu_{(a_n)}$ is defined in (\ref{compensFV}). Then,
since $Y^n_t$ is a purely discontinuous martingale, we have
$B^n_t=(B^{n,i}_t)_{1\leq i\leq p}$, with
$$
B^{n,i}_t=-\int_{[0,t]\times\R^p}
\mu_{(a_n)}(dt,dz_1\ldots dz_p)\,(z_i-h(z_i))
$$
Similarly,
$\tilde C^n_t=(\tilde C^{i,j,n}_t)_{1\leq i,j \leq p}$, with
$$
\tilde C^{i,j,n}_t=\int_{[0,t]\times\R^p}
\mu_{(a_n)}(dt,dz_1\ldots dz_p)h(z_i)h(z_j).$$

Write $C_*(\R^p)$ for the space of all bounded Lipschitz continuous
functions on $\R^p$ that vanish on a neighborhood of $0$. We fix
$g\in C_*(\R^p)$ such that $|g|\leq 1$, and we
choose $\alpha>0$ such that $g(z_1,\ldots,z_p)=0$
if $|z_i|\leq \alpha$ for every $i=1,\ldots,p$. Following the notation
in
\cite{JS}, we set
$$
(g*\mu_{(a_n)})_t =\int_{[0,t]\times\R^p}
\mu_{(a_n)}(dt,dz_1\ldots dz_p)\,g(z_1,\ldots,z_p).$$

 From formula (\ref{compensFV}) we have
\begin{eqnarray}
\label{lapotech0}
B^{n,i}_t&=&\int_0^t ds\,\beta^{n,i}(Y^{n,1}_s,\ldots,Y^{n,p}_s)
\nonumber\\
\tilde C^{n,i,j}_t&=&\int_0^t
ds\,\gamma^{n,i,j}(Y^{n,1}_s,\ldots,Y^{n,p}_s)\\
(g*\mu_{(a_n)})_t &=&\int_0^t
ds\,\varphi^{n}(Y^{n,1}_s,\ldots,Y^{n,p}_s),\nonumber
\end{eqnarray}
where the functions $\beta^{n,i},\gamma^{n,i,j},\varphi^n$
are defined by
\ba
\beta^{n,i}(y_1,\ldots,y_p)
&=&-\int_{\R^p}
R^{(a_n)}(y_1,\ldots,y_p;dz_1,\ldots,dz_p)\,(z_i-h(z_i))\\
\gamma^{n,i,j}(y_1,\ldots,y_p)&=&\int_{\R^p}
R^{(a_n)}(y_1,\ldots,y_p;dz_1,\ldots,dz_p)\,h(z_i)h(z_j)\\
\varphi^n(y_1,\ldots,y_p)
&=&\int_{\R^p}
R^{(a_n)}(y_1,\ldots,y_p;dz_1,\ldots,dz_p)\,g(z_1,\ldots,z_p).
\ea

Similarly, the (modified) characteristics of the semimartingale
$Z$ are
$$(B,\tilde C,\theta)$$
where $\theta$ is as in Proposition \ref{UniCSBP}, and
\begin{eqnarray}
\label{lapotech1}
B^{i}_t&=&\int_0^t ds\,\beta^{i}(Z^{1}_s,\ldots,Z^{p}_s)\nonumber\\
\tilde C^{i,j}_t&=&\int_0^t
ds\,\gamma^{i,j}(Z^{1}_s,\ldots,Z^{p}_s)\\
(g*\theta)_t &=&\int_0^t
ds\,\varphi(Z^{1}_s,\ldots,Z^{p}_s),\nonumber
\end{eqnarray}
where the functions $\beta^{i},\gamma^{i,j},\varphi$
are respectively defined by
\ba
\beta^{i}(y_1,\ldots,y_p)
&=&-\int_{\R^p}
U(y_1,\ldots,y_p;dz_1,\ldots,dz_p)\,(z_i-h(z_i))\\
\gamma^{i,j}(y_1,\ldots,y_p)&=&\int_{\R^p}
U(y_1,\ldots,y_p;dz_1,\ldots,dz_p)\,h(z_i)h(z_j)\\
\varphi(y_1,\ldots,y_p)
&=&\int_{\R^p}
U(y_1,\ldots,y_p;dz_1,\ldots,dz_p)\,g(z_1,\ldots,z_p).
\ea

\begin{lemma}
\label{largepoptech}
For every $(y_1,\ldots,y_p)\in{\cal D}_p$,
\ba
|\beta^{n,i}(y_1,\ldots,y_p)|&\leq&2y_i\int \nu_{(a_n)}(dr)\,r\,
{\bf 1}_{\{r>\delta\}}\\
|\gamma^{n,i,j}(y_1,\ldots,y_p)|&\leq&
(y_i+y_j)\int \nu_{(a_n)}(dr)\,(r\wedge r^2)\\
|\varphi^n(y_1,\ldots,y_p)&\leq& {2\over \alpha}\,y_p\,
\int \nu_{(a_n)}(dr)\,r\,
{\bf 1}_{\{r>\alpha\}}.
\ea
Moreover,
\ba
\lim_{n\to\infty} \beta^{n,i}(y_1,\ldots,y_p)&=&\beta^i(y_1,\ldots,y_p)\\
\lim_{n\to\infty} \gamma^{n,i,j}(y_1,\ldots,y_p)&=&
\gamma^{i,j}(y_1,\ldots,y_p)\\
\lim_{n\to\infty} \varphi^n(y_1,\ldots,y_p)&=&\varphi(y_1,\ldots,y_p),
\ea
uniformly on bounded subsets of ${\cal D}_p$.
\end{lemma}

Let us postpone the proof of the lemma and complete that
of the theorem. The first step is to check the sequence of the
laws of the processes $Y^n$ is tight in the space of probability measures
on $\D(\R_+,\R^p)$. This will follow from Theorem VI.4.18 in \cite{JS}
provided we can check that:
\begin{description}
\item{(i)} We have for every $N>0$ and $\varepsilon>0$,
$$\lim_{b\uparrow\infty} \Big(\limsup_{n\to\infty}
P[\mu_{(a_n)}([0,N]\times\{z\in \R^p:|z|>b\})>\varepsilon]\Big)=0.$$
\item{(ii)} The laws of the processes $B^{n,i},\tilde C^{n,i,j},
g*\mu_{(a_n)}$ are tight in the space of probability measures
on $C(\R_+,\R)$.
\end{description}
To prove (i), set
$$T^n_A=\inf\{t\geq 0:Y^{p,n}_t> A\}$$
for every $A>x_p$. Since $Y^{n,p}$ is a (bounded) nonnegative martingale,
a classical result states that
\be
\label{lptech1}
\P[\sup\{Y^{n,p}_t,t\geq 0\}>A]=\P[T^n_A<\infty]\leq {x_p\over A}.
\ee
 From formulas (\ref{compensFV}) and (\ref{jumpFVbis}), we have on the
event $\{\sup\{Y^{p,n}_t,t\geq 0\}\leq A\}$
$$\mu_{(a_n)}([0,N]\times\{z\in \R^p:|z|>b\})
\leq N\Big(A\,\nu_{(a_n)}(]{b\over p},\infty[)+a_n\,\nu_{(a_n)}
(]{ba_n\over pA},\infty[)\Big)$$
Under Assumption (H), we have
$$\lim_{n\to\infty} a_n\,\nu_{(a_n)}
(]{ba_n\over pA},\infty[) =0$$
and so, on the event $\{\sup\{Y^{p,n}_t,t\geq 0\}\leq A\}$,
$$\limsup_{n\to\infty}
\mu_{(a_n)}([0,N]\times\{z\in \R^p:|z|>b\})
\leq NA\,\pi([{b\over p},\infty[).$$
If we first choose $A$ so that $x_p/A$ is small, and then
$b$ large enough so that $NA\,\pi([{b\over p},\infty[)<\varepsilon$, we see
that the statement in (i) follows from (\ref{lptech1}). Part (ii) is a
straightforward consequence of formulas (\ref{lapotech0}),
the bounds of the first part of Lemma \ref{largepoptech} and (\ref{lptech1})
again. This completes the proof of the tightness of the sequence of the
laws of the processes $Y^n$.

Then, we can assume that, at least along a suitable subsequence,
$Y^n$ converges in distribution towards a limiting process
$Y^\infty=(Y^{\infty,1},\ldots,Y^{\infty,p})$. We claim that
$Y^\infty$ is a semimartingale whose
triplet of (modified) characteristics $(B^\infty,\tilde C^\infty,
\mu_\infty)$ is such that
\begin{eqnarray}
\label{lapotech2}
B^{\infty,i}_t&=&\int_0^t
ds\,\beta^{i}(Y^{\infty,1}_s,\ldots,Y^{\infty,p}_s)\nonumber\\
\tilde C^{\infty,i,j}_t&=&\int_0^t
ds\,\gamma^{i,j}(Y^{\infty,1}_s,\ldots,Y^{\infty,p}_s)\\
(g*\mu_\infty)_t &=&\int_0^t
ds\,\varphi(Y^{\infty,1}_s,\ldots,Y^{\infty,p}_s),\nonumber
\end{eqnarray}
with $\beta^i,\gamma^{i,j},\varphi$ as above. To see this, it is
enough to verify that the 4-tuples $(Y^n,B^n,\tilde C^n,g*\mu_{(a_n)})$
converge in distribution to
$(Y^\infty,B^\infty,\tilde C^\infty,g*\mu_{\infty})$ (see Theorem IX.2.4 in
\cite{JS}). The latter convergence readily follows from the convergence of
$Y^n$ towards $Y^\infty$, formulas  (\ref{lapotech0}) and the second part
of Lemma \ref{largepoptech}.

Finally, knowing the triplet of characteristics
of $Y^\infty$, Theorem II.2.34 in \cite{JS}
shows that $Y^\infty$ is a purely
discontinuous martingale, and the compensator of its jump measure is
$$dt\;U(Y^{\infty,1}_t,\ldots,Y^{\infty,p}_t;dz_1,\ldots,dz_p).$$
By Proposition \ref{UniCSBP}, this implies that
$Y^\infty$ has the same distribution as $Z$, and this completes the
proof of Theorem \ref{largepop}. \QED

\noindent{\bf Proof of Lemma \ref{largepoptech}:} By definition, for
$(y_1,\ldots,y_p)\in {\cal D}_{p}^{a_n}$,
\begin{eqnarray*}
\beta^{n,i}(y_1,\ldots,y_p)
&=&-\int \nu^{(a_n)}(dr)\int_0^{a_n}
du\Big(r({\bf 1}_{\{u\leq y_i\}}-a_n^{-1}y_i)
-h(r({\bf 1}_{\{u\leq y_i\}}-a_n^{-1}y_i))\Big)\\
&=&-\int \nu^{(a_n)}(dr)
\,y_i(r(1-a_n^{-1}y_i)-h(r(1-a_n^{-1}y_i)))\\
&&+\int \nu^{(a_n)}(dr)\,(a_n-y_i)(a_n^{-1}ry_i
+h(-a_n^{-1}ry_i)).
\end{eqnarray*}
Recalling that $h(x)=x$ if $|x|\leq \delta$,
we immediately get the bound
$$\beta^{n,i}(y_1,\ldots,y_p)\leq
2y_i\int \nu^{(a_n)}(dr)\,r{\bf 1}_{\{r>\delta\}}.$$
Furthermore, using the fact that $h$ is Lipschitz
with Lipschitz constant $1$,  we have
\begin{eqnarray*}
&&|\beta^{n,i}(y_1,\ldots,y_p)+y_i\int \nu^{(a_n)}(dr)\,(r-h(r))|\\
&&\qquad\leq 2a_n^{-1}y_i^2 \int \nu^{(a_n)}(dr)\,r{\bf 1}_{\{r>\delta\}}
+y_i\int \nu^{(a_n)}(dr)\,r{\bf 1}_{\{a_n^{-1}ry_i>\delta\}}
\end{eqnarray*}
and it is easy to verify
from Assumption (H) that the right-hand side tends to $0$ as
$n\to\infty$, uniformly when $y_i$ varies over a bounded
subset in $\R_+$. Since Assumption (H) also implies that
$$\lim_{n\to\infty}\int \nu^{(a_n)}(dr)\,(r-h(r))
=\int \pi(dr)\,(r-h(r)),$$
we get the first limit of the lemma.

Consider now, for
$(y_1,\ldots,y_p)\in {\cal D}_{p}^{a_n}$, and $1\leq i\leq j\leq p$,
\begin{eqnarray}
\label{largepoptech1}
\gamma^{n,i,j}(y_1,\ldots,y_p)
&=&\int \nu^{(a_n)}(dr)\int_0^{a_n}
du\,
h(r({\bf 1}_{\{u\leq y_i\}}-a_n^{-1}y_i))
h(r({\bf 1}_{\{u\leq y_j\}}-a_n^{-1}y_j))\Big)\nonumber\\
&=&\int
\nu^{(a_n)}(dr)\,y_i\,h(r(1-a_n^{-1}y_i))h(r(1-a_n^{-1}y_j))\nonumber\\
&&+\int
\nu^{(a_n)}(dr)\,(y_j-y_i)\,h(-a_n^{-1}ry_i)\,h(r(1-a_n^{-1}y_j))\nonumber\\
&&+\int
\nu^{(a_n)}(dr)\,(a_n-y_j)\,h(-a_n^{-1}ry_i)\,h(-a_n^{-1}ry_j).
\end{eqnarray}
Using the bounds $|h|\leq 1$ and $|h(x)|\leq |x|$, we get
$$|\gamma^{n,i,j}(y_1,\ldots,y_p)|
\leq y_j\int \nu^{(a_n)}(dr)(r^2\wedge 1)+y_i\int
\nu^{(a_n)}(dr)\,r(r\wedge 1),$$
which gives the second bound of the lemma. Then, using the
Lipschitz property
of $h$,
\begin{eqnarray*}
&&\Big|\int \nu^{(a_n)}(dr)\,h(r(1-a_n^{-1}y_i))h(r(1-a_n^{-1}y_j))
-\int \nu^{(a_n)}(dr)\,h(r)^2\Big|\\
&&\qquad\leq 2a_n^{-1}y_j\int \nu^{(a_n)}(dr)\,rh(r) \longrightarrow 0
\end{eqnarray*}
as $n\to \infty$. Notice that
$$\lim_{n\to\infty}y_i\int \nu^{(a_n)}(dr)\,h(r)^2
=y_i\int \pi(dr)\,h(r)^2=\gamma^{i,j}(y_1,\ldots,y_p),$$
uniformly when $(y_1,\ldots,y_p)$ varies over a bounded set. To complete
the verification of the second limit in the lemma, we need to
check that the last two terms in the right-hand side of
(\ref{largepoptech1}) tend to $0$ as $n\to\infty$. We have first
$$\int \nu^{(a_n)}(dr)\,h(-a_n^{-1}ry_i)h(r(1-a_n^{-1}y_j))
\leq \int \nu^{(a_n)}(dr)\,ra_n^{-1}y_i\,h(r)\longrightarrow 0$$
as $n\to \infty$. It remains to bound
\begin{eqnarray*}
&&\Big|a_n\int
\nu^{(a_n)}(dr)\,h(-a_n^{-1}ry_i)\,h(-a_n^{-1}ry_j)\Big|\\
&&\quad \leq \int\nu^{(a_n)}(dr)\,ry_i((a_n^{-1}ry_j)\wedge 1)\\
&&\quad \leq y_i\int\nu^{(a_n)}(dr)\,r{\bf 1}_{\{r>A\}}
+y_iy_j a_n^{-1}\int\nu^{(a_n)}(dr)\,r^2{\bf 1}_{\{r\leq A\}}
\end{eqnarray*}
where $A>0$ is arbitrary. If $\eta>0$ is given, we can first choose
$A$ sufficiently large so that
$$\limsup_{n\to\infty} \int\nu^{(a_n)}(dr)\,r{\bf 1}_{\{r>A\}}
\leq \int \pi(dr)\,r{\bf 1}_{\{r\geq A\}}<\eta.$$
On the other hand, we have also
$$\lim_{n\to\infty}a_n^{-1}\int\nu^{(a_n)}(dr)\,r^2{\bf 1}_{\{r\leq A\}}
=0$$
and together with the preceding estimates, this gives the second
limit of the lemma.

Finally, we have
$$\varphi^n(y_1,\ldots,y_p)
=\int\nu^{(a_n)}(dr)\int_0^{a_n} du\, g(
r({\bf 1}_{\{u\leq y_1\}}-a_n^{-1}y_1),\ldots,
r({\bf 1}_{\{u\leq y_p\}}-a_n^{-1}y_p)).$$
Since $|g|\leq 1$ and $g(z_1,\ldots,z_p)=0$ if $\sup|z_i|\leq \alpha$,
we easily get the bound
\begin{eqnarray*}
|\varphi^n(y_1,\ldots,y_p)|&\leq&
y_p\int\nu^{(a_n)}(dr)\,{\bf 1}_{\{r>\alpha\}}
+a_n\int\nu^{(a_n)}(dr)\,{\bf 1}_{\{a_n^{-1}ry_p>\alpha\}}\\
&\leq&y_p\int\nu^{(a_n)}(dr)\,{\bf 1}_{\{r>\alpha\}}
+{y_p\over \alpha}\int\nu^{(a_n)}(dr)\,r{\bf 1}_{\{r>\alpha\}}
\end{eqnarray*}
which gives the third bound of the lemma.
Then, if $M$ denotes a Lipschitz constant for $g$,
\begin{eqnarray*}
&&\Big|\varphi^n(y_1,\ldots,y_p)-\int\nu^{(a_n)}(dr)\int_0^{a_n} du\, g(
r{\bf 1}_{\{u\leq y_1\}},\ldots,
r{\bf 1}_{\{u\leq y_p\}})\Big|\\
&&\quad\leq Mp\int\nu^{(a_n)}(dr)\int_0^{a_n}du
({\bf 1}_{\{u\leq y_p\}}a_n^{-1}ry_p\,{\bf 1}_{\{r>\alpha\}}
+{\bf 1}_{\{u> y_p\}}a_n^{-1}ry_p\,{\bf 1}_{\{a_n^{-1}ry_p>\alpha\}})\\
&&\quad\leq Mp\Big(\int\nu^{(a_n)}(dr)\,r{\bf 1}_{\{r>\alpha\}}\Big)
a_n^{-1}y_p^2+
Mpy_p\int \nu^{(a_n)}(dr)\,r{\bf 1}_{\{a_n^{-1}ry_p>\alpha\}}
\end{eqnarray*}
which tends to $0$ as $n$ tends to $\infty$, uniformly when
$y_p$ varies over a compact subset of $\R_+$. The last
convergence of the lemma now follows from Assumption (H). This completes
the proof. \QED

\section{Hydrodynamic limits for exchangeable coalescents}

The motivation for this section stems from hydrodynamic limit theorems leading
from stochastic coalescents to Smoluchowski's coagulation equation, which we
now summarize.

\subsection{Stochastic coalescents and Smoluchowski's coagulation
equation} Consider a symmetric measurable function $K:]0,\infty[\times
]0,\infty[\to\R_+$ which will be referred to as a coagulation kernel. A
stochastic coalescent with coagulation kernel $K$ can be viewed as a
Markov chain
  in continuous
time $C=(C_t, t\geq0)$ with values in the space of finite integer-valued
measures on
$]0,\infty[$ with the following dynamics. Suppose that the process
starts from some state $\sum_{i=1}^{k}\delta_{x_i}$, where $k\geq 2$ and
$x_i\in]0,\infty[$ for $i=1,\ldots,k$. For $1\leq i < j \leq k$,
let ${\bf e}_{i,j}$ be an exponential variable with
parameter $K(x_i,x_j)$, such that to different pairs correspond
independent variables. The first jump of the process
$C$ occurs at time
$\min_{1\leq i < j \leq k}{\bf e}_{i,j}$, and if this minimum is
reached for the indices $1\leq
\ell < m
\leq k$ (i.e. $\ell$ and $k$ are the indices such that
$\min_{1\leq i < j \leq k}{\bf e}_{i,j}={\bf e}_{\ell,m}$),
then the state after the jump is
$$\delta_{x_\ell+x_m}+\sum_{i\neq \ell, m}\delta_{x_i}\,.$$
In other words, a stochastic coalescent with coagulation kernel $K$
is a finite particle system in $]0,\infty[$ such that
each pair of particles $(x_i,x_j)$ in the system merges at rate
$K(x_i,x_j)$, independently of the other pairs.

Now consider a sequence $(\tilde C^{(n)}_t, t\geq0)_{n\in\N}$ of stochastic
coalescents with coagulation kernel $K$ and set
$C^{(n)}_t=n^{-1}\tilde C^{(n)}_{t/n}$ for $t\geq0$.
Suppose that the sequence of initial states $C^{(n)}_0$
converges in probability in ${\cal M}_{\rm R}$ to a Radon measure $\mu_0$. Then
under some technical assumptions on the coagulation kernel $K$ (see e.g.
Norris
\cite{Norris}), the sequence  $(C^{(n)}_t, t\geq0)$ converges in
probability on the space of c\`
adl\` ag trajectories with values in ${\cal M}_{\rm R}$ towards a
deterministic limit
$(\mu_t, t\geq0)$.
Moreover this limit is characterized as the solution to
Smoluchowski's coagulation
equation
\begin{equation}\label{eqsmolu}
{d \langle \mu_t,f \rangle \over dt}
\,=\,{1\over
2}\int_{]0,\infty[^2}
\left(f(x+y)-f(x)-f(y)\right)K(x,y)\mu_t(dx)\mu_t(dy)\,,
\end{equation}
where $f:]0,\infty[\to\R$ denotes a generic continuous function with
compact support.

\subsection{Hydrodynamic limits}

Let  $\nu$ denote a
$\sigma$-finite measure on
$]0,1]$ such that $\int_{]0,1]}r^2\nu(dr)<\infty$, and
write $\Lambda(dr)=r^2\nu(dr)$, which is thus a
finite measure on $]0,1]$. The so-called $\Lambda$-coalescent
(or coalescent
with multiple collisions, see \cite{Pi}) is a Markov process
$(\Pi_t,t\geq 0)$
taking values in the set of all partitions of $\N$.
Unless otherwise specified, we assume that $\Pi_0$ is the partition of
$\N$ into singletons. For every $t\geq 0$, write $D_t$ for the sequence
of asymptotic frequencies of the blocks of $\Pi_t$, ranked
in nonincreasing order (if the number $k$ of blocks is finite, then the terms
of index greater than $k$ in the sequence are all equal to $0$).
Then (\cite{Pi}, section 2.2) the process $(D_t,t\geq 0)$
is a time-homogeneous Markov
process with values in the space ${\cal S}^{\downarrow}_1$
of nonincreasing numerical sequences ${\bf s}=(s_1,\ldots)$
with $\sum_{i=1}^{\infty}s_i\leq 1$.

The following connection with generalized Fleming-Viot
processes can be found in \cite{BeLG1}. Let $(F_t,t\geq 0)$ be the generalized
Fleming-Viot process associated with $\nu$, and for every $t\geq 0$, let
$J_t$ be the sequence of sizes of jumps of the mapping
$x\to F_t(x)$, ranked again in nonincreasing order, and with the
same convention if there are finitely many jumps. Then, for each fixed
$t\geq 0$, $J_t$ and $D_t$ have the same distribution (Theorem 1 in
\cite{BeLG1} indeed gives a deeper connection, which has been
briefly described in Section 1).

For each $a>0$, let $\tilde \nu^{(a)}$, $\nu^{(a)}$,
$\tilde F^{(a)}$ and $F^{(a)}$ be as in Section 3. Denote by
$\tilde\mu^{(a)}_t$ the point measure whose atoms are given by the
jump sizes of the increasing process $x\to \tilde F^{(a)}_t(x)$:
$$\tilde\mu^{(a)}_t=\sum_{\{x\in]0,1]:\tilde F^{(a)}_t(x)-\tilde
F^{(a)}_t(x-)>0\}}
\delta_{\tilde F^{(a)}_t(x)-\tilde F^{(a)}_t(x-)}.$$
Fom the preceding observations, the atoms of $\tilde \mu^{(a)}_t$
also correspond
to the sizes of the blocks in a $\Lambda$-coalescent at time $t$, for
$\Lambda(dr)=r^2\tilde\nu^{(a)}_t(dr)$.
We then consider the rescaled version
$\mu^{(a)}_t$, given as the image of $a^{-1}\tilde \mu^{(a)}_{at}$ under the
dilation $r\to ar$. Equivalently, $\mu^{(a)}$ is $a^{-1}$ times the
sum of the Dirac
point masses at the jump sizes of the mapping $x\to F^{(a)}_t(x)$.

\begin{theorem}\label{T3} Suppose that {\rm(H)} holds and let
$(Z(t,x);t,x\geq 0)$ be the flow of continuous-state branching
processes associated with $\pi$. Then for every
$t\geq 0$, $\mu^{(a)}_t$ converges to the L\'evy measure $\lambda_t$
of the subordinator $Z(t,\cdot)$
as $a\to\infty$ in probability in ${\cal
M}_{\rm R}$.
\end{theorem}

Theorem \ref{T3} is an immediate consequence of Theorem \ref{largepop}
and the following lemma.

\begin{lemma}\label{L2} Let
$\sigma=(\sigma_t, t\geq0)$ be a subordinator with L\'evy
measure $\lambda$.
For each $a>0$, let $X^{(a)}=(X^{(a)}_t, 0\leq t
\leq a)$ be an increasing c\` adl\` ag process with exchangeable
increments, with $X^{(a)}_0=0$ and $X^{(a)}_a=a$ a.s.
Suppose that $X^{(a)}$ converges to $\sigma$ as
$a\to\infty$ in the sense of finite-dimensional distributions.
Then the random point measure
$$a^{-1}\sum_{0<t<a}\delta_{\Delta X^{(a)}_t}$$
converges to $\lambda$ in probability in ${\cal
M}_{\rm R}$ as $a\to\infty$. \end{lemma}

\proof Pick some nonnegative continuous function $f:]0,\infty[\to\R$ with
compact support and write
$$c:=\int_{]0,\infty[}f(x)\lambda (dx)\,.$$

By the L\'evy-It\^o decomposition for subordinators,
the random point measure $\sum_{\Delta\sigma_t>0}\delta_{(t,
\Delta\sigma_t)}$ on $\R_+\times ]0,\infty[$ is Poisson with intensity
$dt\otimes
\lambda(dx)$.
Let $\rho>0$.
The law of large numbers ensures the existence of a real number
$a_{\rho}>0$ such that
\be
\label{L2tech1}
\E\Big[\Big|a_\rho^{-1}\sum_{0<t<
a_\rho}f(\Delta\sigma_t) - c\Big|\Big]< \rho\,.
\ee

Then consider for $a>a_\rho$ the
bridges  with exchangeable increments, bounded variation
and no negative jumps on the time interval $[0, a_\rho]$, defined by
$$B^{(a)}_t:=X^{(a)}_t-t a_{\rho}^{-1}X^{(a)}_{a_\rho}\quad
,\quad B_t=\sigma_t-ta_{\rho}^{-1}\sigma_{a_\rho}\,,
\qquad t\in[0,a_\rho]\,.$$
Our assumptions entail that $B^{(a)}$ converges in the sense
of finite dimensional distributions to $B$, so according to Kallenberg
\cite{Ka}, the random measure
$$\sum_{0<t<a_\rho}\delta_{\Delta
B^{(a)}_t} = \sum_{0<t<a_\rho}\delta_{\Delta X^{(a)}_t}$$
converges in law on ${\cal M}_{\rm R}$ towards
$$\sum_{0<t<a_\rho}\delta_{\Delta B_t}\,
=\,\sum_{0<t<a_\rho}\delta_{\Delta\sigma_t}\,,$$
and in particular,
when $a\to\infty$,
\be
\label{L2tech2}
a_\rho^{-1}\sum_{0<t<
a_\rho}f(\Delta X^{(a)}_t)\
\build{\la}_{}^{\rm (d)} \ a_\rho^{-1}\sum_{0<t<
a_\rho}f(\Delta\sigma_t)\,.
\ee

Let us check that the variables
$$\Big|\sum_{0<t<
a_\rho}f(\Delta X^{(a)}_t)\Big|\,,\qquad a\in[a_\rho,\infty[
$$
are uniformly
integrable. Let $[u,v]$ be a compact subinterval of $]0,\infty[$
such that the support of $f$ is contained in $[u,v]$. Denote by
$N^a_{[u,v]}$ the number of jumps of the process $X^{(a)}$
with size in $[u,v]$. By classical results about processes
with exchangeable increments, conditionally on
$N^a_{[u,v]}=n$, the number
$$N^{(a,a_\rho)}_{[u,v]}:=\sum_{0<t<
a_\rho} {\bf 1}_{[u,v]}(\Delta X^{(a)}_t)$$
has a binomial ${\cal B}(n,{a_\rho\over a})$ distribution. Notice that
$N^a_{[u,v]}\leq {a\over u}$ since $X^{(a)}_a=a$.
We see that $N^{(a,a_\rho)}_{[u,v]}$ is bounded above in
distribution by a binomial ${\cal B}([{a\over u}],{a_\rho\over a})$
distribution, and the desired uniform integrability readily follows.

It then follows from (\ref{L2tech1}) and (\ref{L2tech2}) that
$$\lim_{a\to\infty}
\E\Big[\Big|a_\rho^{-1}\sum_{0<t<
a_\rho}f(\Delta X^{(a)}_t)-c\Big|\Big]
\,=\,\E\Big[\Big|a_\rho^{-1}\sum_{0<t<
a_\rho}f(\Delta\sigma_t) - c\Big|\Big]\leq \rho\,.
$$
Moreover, an easy exchangeability argument shows that we have also
$$\limsup_{a\to\infty}
\E\Big[\Big|a^{-1}\sum_{0<t<
a}f(\Delta X^{(a)}_t)-c\Big|\Big]\leq \rho\,.$$
Since $\rho$ may be taken arbitrarily small, we have
thus shown that
$$\lim_{a\to\infty}
a^{-1}\sum_{t\leq
a}f(\Delta X^{(a)}_t)\,=\,\int_{]0,\infty[}f(x)\lambda(dx)\,,$$
in $L^1$ for every continuous function $f$ with compact support.
The conclusion now follows by a standard argument. \QED

We will now show that the
family $(\lambda_t, t>0)$ of L\'evy measures, which appears in Theorem
\ref{T3}, solves a certain coagulation equation with multiple collisions.
To this end, we introduce the following additional assumption, which
also plays a key role
in the study of the genealogical structure of continuous-state
branching processes
(see e.g. \cite{DuLG}).

\noindent{\bf Assumption} (E) {\it The $\Psi$-CSBP becomes extinct
almost surely.}

Equivalently, this assumption holds iff $\P[Z(t,x)=0]>0$ for every
$t>0$ and $x\geq 0$.
By solving (\ref{eqsg}), it is easy to verify that Assumption (E)
is equivalent to
\be
\label{extinctcond}
\int_1^\infty {du\over \Psi(u)}<\infty.
\ee
In particular, Assumption (E) holds in the so-called stable case
$\Psi(u)=u^\gamma$, $\gamma\in]1,2[$, that will be considered in
Section 5 below.

 From (\ref{driftLevy}), we see that under Assumption (E) we have
$d_t=0$, and the total mass $\lambda_t(]0,\infty[)=-\log \P[Z(t,1)=0]$ is
finite for every $t>0$. Moreover the function $t\to\lambda_t(]0,\infty[)$
is nonincreasing.

We denote by $C_{\bullet}(\R_+)$ the space of all
bounded continuous functions $f$ on $\R_+$ such that $f(0)=0$
and $f(x)$ has a limit as $x\to+\infty$. The space $C_{\bullet}(\R_+)$
is equipped with the uniform norm, which is denoted by $\|f\|$. For
every integer $k\geq 2$ and $q>0$, we denote
by $\Psi^{(k)}(q)$ the $k$-th derivative of $\Psi$ at $q$. It is
immediately checked
that
\be
\label{derivpsi}
\Psi^{(k)}(q)= (-1)^k\int \pi(dr)\,r^k\e^{-qr}.
\ee
Obviously, $(-1)^k\Psi^{(k)}(q)\geq 0$ for every $k\geq 2$ and $q>0$.

\begin{proposition}
\label{Smolu}
Under Assumption {\rm (E)}, for every $f\in C_{\bullet}(\R_+)$,
the function $t\to \langle \lambda_t,f\rangle$ solves the
equation
\be
\label{Smolueq0}
{d \langle\lambda_t,f\rangle\over dt}=
\sum_{k=2}^\infty {(-1)^k\Psi^{(k)}(\langle \lambda_t,1\rangle)\over k!}
\int_{]0,\infty[^k}
(f(x_1+\cdots+x_k)-(f(x_1)+\cdots+f(x_k)))\,\lambda_t(dx_1)\ldots
\lambda_t(dx_k)
\ee
where the series in the right-hand side converges absolutely.
\end{proposition}

It is interesting to observe that (\ref{Smolueq0}) also holds when
$\Psi(q)=cq^2$ for some
constant $c>0$. Take $\Psi(q)={1\over 2}u^2$ for definiteness
(then the $\Psi$-CSBP is the
classical Feller diffusion)  in such a way that
(\ref{Smolueq0}) exactly reduces to (\ref{eqsmolu}) with $K\equiv 1$.
Then $u_t(q)=2q\,(2+qt)^{-1}$,
and it follows that
\be
\label{Fellerdens}
\lambda_t(dx)={4\over t^2}\exp(-{2x\over t})\,dx
\ee
so that the density of $\lambda_t$ is the classical solution,
arising from infinitesimally small initial clusters, of the
Smoluchovski equation
(\ref{eqsmolu}) in the case $K\equiv 1$ (cf Section 2.2 of
\cite{Al}).

We can rewrite equation (\ref{Smolueq0}) in a somewhat more synthetic
way by introducing the following notation.
If $\mu$ is a measure on  $]0,\infty[$
such that
$\int_{]0,\infty[} (1\wedge x)\mu(dx)<\infty$, we write $\mu^{\oplus}$ for
the distribution on $[0,\infty[$ of the sum of the atoms of a Poisson
random measure on $]0,\infty[$ with intensity $\mu$. Note that
$\mu^{\oplus}$ is a probability measure and that, by Campbell's formula,
\be
\label{Camp}
\int_{[0,\infty[}\e^{-qx}\mu^{\oplus}(dx)\,=\,
\exp\Big\{-\int_{]0,\infty[}(1-\e^{-qx})\mu(dx)\Big\}\,,\qquad
q\geq0.
\ee
As we will see in the proof below, (\ref{Smolueq0}) follows from the equation
\be
\label{Smolueq}
{d \langle\lambda_t,f\rangle\over dt}\,=\,\int_{]0,\infty[}\pi(da)
\Big(\langle(a\lambda_t)^{\oplus},f\rangle-\langle a\lambda_t,f\rangle\Big).
\ee
Informally, we may think of $\lambda_t(dx)$ as the density at time $t$ of
particles with size $x$ in some infinite system of particles. The
right-hand side in (\ref{Smolueq}) can be interpreted by saying
at rate $\pi(da)$, a `quantity' $a$ of particles coagulates at time $t$.
More precisely, this `quantity' is sampled in a Poissonian way, viewing at
$a\lambda_t$ as an intensity measure for the sampling (so, loosely speaking,
the particles involved into the coagulation are sampled uniformly at
random amongst the particles present at time $t$).

As the proof below will show, (\ref{Smolueq}) still holds without
Assumption (E)
at least for functions $f$ of the type
$f(x)=1-\exp(-q x)$, provided that $d_t=0$ for every $t>0$
(recall from Silverstein \cite{Silv} that the latter holds whenever
$\int_{]0,1[}r\pi(dr)=\infty$).
In that case however, the measures $\lambda_t$ may be infinite, and then
coagulations involve infinitely many components, so
that one cannot write an equation of the form (\ref{Smolueq0}).

\proof We first prove (\ref{Smolueq}). For $q>0$, let $f_{(q)}\in
C_{\bullet}(\R_+)$
be defined by $f_{(q)}(x)=1-\e^{-q x}$.
By (\ref{Camp}) and (\ref{driftLevy}),
$$\langle(a\lambda_t)^{\oplus},f_{(q)}\rangle
=1-\exp\Big(-a\int\lambda_t(dr)(1-\e^{-qr})\Big)=1
-\exp(-a u_t(q)).$$
On the other hand, by (\ref{driftLevy}) again,
$$\langle \lambda_t,f_{(q)}\rangle=u_t(q).$$
Thus when $f=f_{(q)}$ the right-hand side of (\ref{Smolueq}) makes sense
and is equal to
$$\int_{]0,\infty[}\pi(da)
\Big(1-\exp(-a u_t(q))-au_t(q)\Big)=-\Psi(u_t(q)).$$
Therefore (\ref{Smolueq}) reduces to (\ref{eqsg}) in that case. Note that
we have not used Assumption (E) at this stage (except for the fact that
$d_t=0$ for every $t>0$).

Denote by ${\cal H}$ the subspace of $C_{\bullet}(\R_+)$ that
consists of linear combinations of the functions $f_{(q)}$. Then
${\cal H}$ is dense in $C_{\bullet}(\R_+)$. Obviously,
for every $f\in {\cal H}$,
(\ref{Smolueq}) holds, and the right-hand side
of (\ref{Smolueq}) is a continuous function of $t\in]0,\infty[$.
Fix $f\in C_{\bullet}(\R_+)$ and a sequence $(f_n)_{n\geq 1}$
in ${\cal H}$ that converges to $f$. If we also fix $0<\varepsilon<t$,
we have for every $n\geq 1$,
\be
\label{Smolutech1}
\langle\lambda_t,f_n\rangle=\langle\lambda_\varepsilon,f_n\rangle
+\int_\varepsilon^t ds\int \pi(da)
\Big(\langle(a\lambda_s)^{\oplus},f_n\rangle-\langle
a\lambda_s,f_n\rangle\Big).
\ee
Plainly, for every $s>0$,
$$\langle\lambda_s,f_n\rangle\build{\longrightarrow}_{n\to\infty}^{}
\langle\lambda_s,f\rangle
\quad\hbox{and}\quad
\langle(a\lambda_s)^\oplus,f_n\rangle\build{\longrightarrow}
_{n\to\infty}^{}\langle(a\lambda_s)^\oplus,f\rangle.$$
We claim that there exists a constant $C_\varepsilon$ such that, for every
$s\geq\varepsilon$ and $n\geq 1$, and every $h\in C_\bullet(\R_+)$,
\be
\label{Smolutech2}
|\langle(a\lambda_s)^\oplus,h\rangle - \langle
a\lambda_s,h\rangle|\leq C_\varepsilon(a^2\wedge a)\,\|h\|.
\ee
As the quantities $\langle \lambda_s,1\rangle$,
$s\in[\varepsilon,\infty[$ are bounded above, it is clearly
enough to consider $a\leq 1$. Since $h(0)=0$, the definition of
$(a\lambda_s)^\oplus$ immediately gives
$$\langle(a\lambda_s)^\oplus,h\rangle=\e^{-a\langle\lambda_s,1\rangle}
\,a\langle\lambda_s,h\rangle +O(a^2\|h\|)$$
where the remainder $O(a^2\|h\|)$, which corresponds to the event
that a Poisson measure
with intensity $a\lambda_s$ has at least two atoms, is uniform in
$h\in C_\bullet(\R_+)$ and $s\geq \varepsilon$.
The estimate (\ref{Smolutech2}) follows.

Using (\ref{Smolutech2}) and dominated convergence, we get
\be
\label{Smolutech3}
\lim_{n\to\infty}\int
\pi(da)\Big(\langle(a\lambda_s)^{\oplus},f_n\rangle-\langle
a\lambda_s,f_n\rangle\Big)
=\int \pi(da)
\Big(\langle(a\lambda_s)^{\oplus},f\rangle-\langle a\lambda_s,f\rangle\Big)
\ee
uniformly in $s\in[\varepsilon,\infty[$, and the right-hand side of
(\ref{Smolutech3})
is a continuous function of $s$. Equation (\ref{Smolueq}) in the
general case follows
by passing to the limit $n\to\infty$ in (\ref{Smolutech1}).

Then, to derive (\ref{Smolueq}) from (\ref{Smolueq0}), we write
\ba
&&\int_{]0,\infty[}\pi(da)
\Big(\langle(a\lambda_t)^{\oplus},f\rangle-\langle a\lambda_t,f\rangle\Big)\\
&&=\int \pi(da)\Big(\Big(\sum_{k=1}^\infty
{a^k\over k!}\,\e^{-a\langle\lambda_t,1\rangle}
\int
f(x_1+\cdots+x_k)\,\lambda_t(dx_1)\ldots\lambda_t(dx_k)\Big)-a\langle
\lambda_t,f\rangle\Big)\\
&&=\int \pi(da)\sum_{k=1}^\infty
{a^k\over k!}\,\e^{-a\langle\lambda_t,1\rangle}
\int
(f(x_1+\cdots+x_k)-(f(x_1)+\cdots+f(x_k)))\,\lambda_t(dx_1)\ldots
\lambda_t(dx_k).
\ea
Notice that the term $k=1$ in the last series vanishes. Moreover,
bounding the other terms
by their absolute value gives a convergent series, whose sum is integrable
with respect to $\pi(da)$. Hence we may interchange the sum and the
integral with respect
to $\pi(da)$, and we get the statement of the proposition from (\ref{Smolueq}).
\QED

\noindent{\bf Remark.} To conclude this section, let us observe that
Assumption (E) is closely related to the property for a $\Lambda$-coalescent to
come down from infinity (cf Pitman \cite{Pi} and Schweinsberg \cite{Sch}).
Let  $\nu$ denote a
$\sigma$-finite measure on
$]0,1]$ such that $\int_{]0,1]}r^2\nu(dr)<\infty$, and let
$\Lambda(dx)=x^2\nu(dx)$.
Let $\Psi$ be given by (\ref{mechbranch}) with $\pi=\nu$ (and $\beta=0$).
Then the $\Lambda$-coalescent comes down from infinity if
and only if the $\Psi$-CSBP becomes extinct almost surely.
To see this, recall from Schweinsberg \cite{Sch} that a
necessary and sufficient condition
for the $\Lambda$-coalescent to come down from infinity is
\be
\label{schcond}
\sum_{b=2}^\infty \Big(\sum_{k=2}^b (k-1)\left(
\begin{array}{ll}
b\\
k
\end{array}
\right)\int r^k(1-r)^{b-k}\nu(dr)\Big)^{-1}<\infty.
\ee
Using the binomial formula, we can rewrite this condition as
$$\sum_{b=2}^\infty \Big(
\int (br-1+(1-r)^b)\,\nu(dr)\Big)^{-1}<\infty,
$$
or equivalently, if we put $\Phi(q)=\int (qr-1+(1-r)^q)\,\nu(dr)$ for
every real $q\geq 1$,
\be
\label{schcond2}
\int_2^\infty {dq\over \Phi(q)}<\infty.
\ee
(note that the function $\Phi$ is nondecreasing on $[1,\infty[$).
Simple estimates give the existence of a constant $c\in]0,1[$ such
that, for every
$q\geq 2$,
$$c\Psi(q)\leq \Phi(q)\leq \Psi(q).$$
It follows that (\ref{schcond}) and (\ref{extinctcond}) are equivalent. In the
spirit of the present work, it would be interesting to give a direct
probabilistic proof
of the equivalence between  the property for a
$\Lambda$-coalescent to come down from infinity and Assumption (E)
for the associated branching process.

\section{Small time behavior of flows and coalescents}

In this section, we fix a measure $\nu$ on $]0,1]$
such that $\int r^2\nu(dr)<\infty$ and we consider the
associated generalized Fleming-Viot process $(F_t, t\geq0)$.

 From now on until the end of the section, we
make the following assumption on $\nu$.

\noindent{\bf Assumption} (A). {\it The function
$\nu([\varepsilon,1])$ is regularly varying
with index $-\gamma$ as $\varepsilon\to 0$, for some $\gamma\in]1,2[$.}

As a consequence, there exists a function $L(\varepsilon)$,
$\varepsilon\in]0,1]$ that is
slowly varying as $\varepsilon\to 0$, such that, for every
$\varepsilon\in]0,1]$,
$$\nu([\varepsilon,1])=\varepsilon^{-\gamma}L(\varepsilon).$$

Fix $\varepsilon_0>0$ such that $\nu([\varepsilon_0,1])>0$. For
$\varepsilon\in]0,\varepsilon_0]$ we have
$L(\varepsilon)>0$ and so we can set
$$F^\varepsilon_t(x)={1\over
\varepsilon}\,F_{L(\varepsilon)^{-1}\varepsilon^{\gamma-1}t}(\varepsilon x)
$$
for $0\leq x\leq \varepsilon^{-1}$ and $t\geq 0$. We also let
$\nu_\varepsilon$ be the measure on $[0,\varepsilon^{-1}]$ defined by
$$\int
\nu_\varepsilon(dr)\,\varphi(r)=L(\varepsilon)^{-1}\varepsilon^\gamma
\int\nu(dr)\,\varphi({r\over \varepsilon}).$$
A simple scaling transformation shows that for every
$(x_1,\ldots,x_p)\in{\cal D}_p^{1/\varepsilon}$,
$(F^\varepsilon_t(x_1),\ldots,F^\varepsilon_t(x_p))$ is
a purely discontinuous martingale, with values in ${\cal
D}^{1/\varepsilon}_p$, and the compensator of its jump measure
is
$$dt\,R_\varepsilon(F^\varepsilon_t(x_1),\ldots,F^\varepsilon_t(x_p);
dz_1,\ldots,dz_p)$$
where
\begin{eqnarray}
\label{jumpFV}
&&\int
R_\varepsilon(y_1,\ldots,y_p;dz_1,\ldots,dz_p)\,\varphi(z_1,\ldots,z_p)
\nonumber\\
&&\quad =\int\nu_\varepsilon(dr)\int_0^{1/\varepsilon} du\,\varphi(r({\bf
1}_{\{u\leq y_1\}}-\varepsilon y_1),\ldots,r({\bf 1}_{\{u\leq
y_p\}}-\varepsilon y_p)).
\end{eqnarray}

Let $\pi_\gamma$ be the measure on $]0,\infty[$ such that
$\pi_\gamma(]a,\infty[)=a^{-\gamma}$ for every $a>0$, and let
$$\Psi_\gamma(q)=\int
\pi_\gamma(dr)\,(\e^{-qr}-1+qr)={\Gamma(2-\gamma)\over \gamma-1}
\,q^\gamma.$$
We let $(Z(t,x),t\geq 0,x\geq 0)$
be the flow of continuous-state branching processes constructed
in Section 2, with $\Psi=\Psi_\gamma$.

\begin{theorem}
\label{smalltime}
Under Assumption {\rm (A)}, for every $(x_1,\ldots,x_p)\in{\cal D}_p$,
$$((F^\varepsilon_t(x_1),\ldots,F^\varepsilon_t(x_p));t\geq 0)
\build{\la}_{\varepsilon\to 0}^{\rm(d)} ((Z(t,x_1),\ldots,Z(t,x_p));
t\geq 0)$$
in the Skorokhod space $\D(\R_+,\R^p)$.
\end{theorem}

\proof
This is a simple consequence of Theorem \ref{largepop}, or rather
of its proof. Indeed,
we immediately see that the kernel
$R_\varepsilon(y_1,\ldots,y_p;dz_1,\ldots,dz_p)$ coincides with
$R^{(1/\varepsilon)}(y_1,\ldots,y_p;dz_1,\ldots,dz_p)$ defined in
(\ref{jumpFVbis}), provided we take
$\nu^{(1/\varepsilon)}=\nu_\varepsilon$. From the observation at the
beginning of the proof of
Theorem \ref{largepop}, we see that Theorem \ref{smalltime} will follow
if we can check  that Assumption (H) holds in the present setting, that is
if
\be
\label{keysmall}
\lim_{\varepsilon\to 0} (r\wedge r^2)\,\nu_\varepsilon(dr)
=(r\wedge r^2)\,\pi_\gamma(dr)
\ee
in the sense of weak convergence in ${\cal M}_{\rm F}$.

In order to prove (\ref{keysmall}), first note that when $\varepsilon\to 0+$,
$$\int_{[0,\varepsilon]}x^2\nu(dx)=2\int_{0}^{\varepsilon}
y\nu([y,1])dy-\varepsilon^2\nu([\varepsilon,1])\sim
{\gamma\over 2-\gamma}\varepsilon^{2-\gamma}L(\varepsilon),$$
where the equivalence follows from Assumption (A)
and a classical property of integrals of regularly varying functions.
We immediately deduce that
\be
\label{keysmallvague}
\lim_{\varepsilon\to 0}  r^2\,\nu_\varepsilon(dr)
=r^2\,\pi_\gamma(dr)
\ee
in the sense of {\it vague} convergence in the
space of Radon measures on $[0,\infty[$.

Next, note that
\be
\label{smalltech3}
\int \nu_\varepsilon(dr)\,(r-r\wedge a)
=\int_a^\infty dr\,\nu_\varepsilon([r,\infty[)
=\varepsilon^\gamma L(\varepsilon)^{-1}
\int_a^\infty dr\,\nu([r\varepsilon,1])
\build{\la}_{a\to\infty}^{} 0
\ee
uniformly in $\varepsilon\in]0,\varepsilon_0]$.
 From (\ref{keysmallvague}) and (\ref{smalltech3})
the family $((r\wedge r^2)\,\nu_\varepsilon(dr), 0<\varepsilon <1)$
is tight for the weak topology in ${\cal M}_{\rm F}$.
Together, with (\ref{keysmallvague}), this
establishes the weak convergence (\ref{keysmall}). \QED

\noindent{\bf Remark.} Suppose that $(F_t,t\geq 0)$ is the flow of
bridges associated
with the Kingman coalescent, corresponding to $\Lambda=\delta_0$ in
our notation
(cf Section 4 in \cite{BeLG2}). If we fix $(y_1,\ldots,y_p)\in{\cal
D}^1_p$, the
process $(F_t(y_1),\ldots,F_t(y_p))$ is a diffusion process in ${\cal D}^1_p$
with generator
$${\cal A}g(x)= {1\over 2} \sum_{i,j=1}^p x_{i\wedge j}(1-x_{i\vee
j}){\partial^2 g\over \partial
x_i\partial x_j}(x)$$
(see Theorem 3 in \cite{BeLG2}). Putting $F^\varepsilon_t(x)={1\over
\varepsilon}F_{\varepsilon t}(\varepsilon x)$,
it is a simple matter to verify that our Theorem \ref{smalltime}
still holds in that setting,
provided we let $(Z(t,x),t\geq 0,x\geq 0)$ be the flow associated
with the Feller diffusion
($\Psi(q)={1\over 2}q^2$). Indeed, if we specialize to the case $p=1$
and if we let $(B_t,t\geq 0)$ be a standard linear Brownian motion,
this is just saying
that, for the Fisher-Wright diffusion $(X_t(x),t\geq 0)$ solving
$$dX_t=\sqrt{X_t(1-X_t)}\,dB_t\quad,\quad X_0=x,$$
the rescaled processes $X^\varepsilon_t:={1\over
\varepsilon}X_{\varepsilon t}(\varepsilon x)$
converge in distribution as $\varepsilon\to 0$ towards the Feller
diffusion $Y_t(x)$ solving
$$dY_t=\sqrt{Y_t}\,dB_t\quad,\quad Y_0=x.$$

We will now use Theorem \ref{smalltime} to derive precise information
on the sizes of blocks in a $\Lambda$-coalescent (for $\Lambda(dr)=r^2\nu(dr)$)
in small time.
As previously, we denote by $\lambda_1(dr)$ the L\'evy measure of the
subordinator $(Z(1,x),x\geq 0)$. We have for every $q\geq 0$
$$\exp-x\int (1-\e^{-qr})\lambda_1(dr)=\E(\exp-qZ(1,x))=\exp-xu_1(q)$$
and the function $u_1(q)$ can be calculated from equation
(\ref{eqsg}), with $\Psi=\Psi_\gamma$. It
follows that
$$\int
(1-\e^{-qr})\lambda_1(dr)=(\Gamma(2-\gamma)+q^{1-\gamma})^{1/(1-\gamma)}$$
and in particular, the total mass of $\lambda_1$ is
$$(\Gamma(2-\gamma))^{1/(1-\gamma)}.$$
We will need the fact that $\lambda_1$ has no atoms. An easy way to
derive this property
is to argue by contradiction as follows. Suppose that $a>0$ is an
atom of $\lambda_1$. From the L\'evy-Khintchin
decomposition of $Z(t,x)$ (see the discussion after
(\ref{driftLevy})), it follows that
$a$ is also an atom of the distribution of $Z(1,x)$, for every $x>0$.
By a simple scaling argument,
for every $s>0$, the image of $\lambda_1(dr)$ under the mapping $r\to
s^{1/(\gamma-1)}r$
is $s^{1/(\gamma-1)}\lambda_s(dr)$. Therefore, for every $s\in]0,1[$,
$s^{1/(\gamma-1)}a$
is also an atom of $\lambda_s$, hence of the distribution of $Z(s,x)$
for every $x>0$.
However, applying the Markov property to the process $(Z(t,1))_{t\geq
0}$ at time $1-s$,
this would imply that for every $s\in]0,1[$, $s^{1/(\gamma-1)}a$ is
an atom of the
distribution of $Z(1,1)$, which is absurd.

We set $g(\varepsilon)=L(\varepsilon)^{-1}\varepsilon^{\gamma-1}$
for every $\varepsilon\in]0,\varepsilon_0]$.

\begin{theorem}
\label{blocks}
Assume that {\rm(A)} holds and let $\Lambda(dr)=r^2\nu(dr)$.
For every $t\geq 0$ and $r\in[0,\infty]$, denote by $N_t(]0,r[)$
the number of blocks at time $t$ with frequencies less than $r$ in a
$\Lambda$-coalescent started from the partition of $\N$ in singletons.
Then,
$$\sup_{x\in]0,\infty[}\Big|\varepsilon
N_{g(\varepsilon)}(]0,\varepsilon x[) -\lambda_1 (]0,x[)\Big|
\build{\la}_{\varepsilon\to 0}^{} 0$$
in probability.
\end{theorem}

Again Theorem \ref{blocks} is a generalization of classical results
for the Kingman coalescent. In that case, one has
$$\sup_{x\in]0,\infty[}\Big|\varepsilon
N_{\varepsilon}(]0,\varepsilon x[) - 2(1-2\e^{-2 x})\Big|
\build{\la}_{\varepsilon\to 0}^{} 0$$
almost surely (cf Section 4.2 of \cite{Al}). This is consistent
with Theorem \ref{blocks} since in the case $\Psi(q)={1\over 2}q^2$,
(\ref{Fellerdens}) shows that
$$2(1-2\e^{-2 x})=\int_0^x 4\e^{-2x}\,dx=\lambda_1(]0,x[).$$

\proof By the results of \cite{BeLG1} recalled at the beginning of
subsection 4.2, we know that,
for each $t\geq 0$ the collection
$(N_t(]0,r[),r\geq 0)$ has the same distribution as
$$\Big(\sum_{0<u<1} {\bf 1}_{\{0<F_t(u)-F_t(u-)\leq r\}}\,,\,r\geq 0\Big),$$
where $(F_t,t\geq 0)$ is the generalized Fleming-Viot process
associated with $\nu$.
It then follows from our definitions that
$$\Big(\varepsilon\,N_{g(\varepsilon)}(]0,x\varepsilon[),x\geq 0\Big)
\build{=}_{}^{\rm(d)}
\Big(\varepsilon\sum_{0<u<1/\varepsilon}
{\bf 1}_{\{0<F^\varepsilon_1(u)-F^\varepsilon_1(u-)\leq
x\}}\,,\,x\geq 0\Big).$$
By combining Theorem \ref{smalltime} and Lemma \ref{L2}, we get
that
\be
\label{blocks00}
\varepsilon\sum_{0<u<1/\varepsilon}\delta_{F^\varepsilon_t(u)-
F^\varepsilon_t(u-)}
\build{\la}_{\varepsilon\to 0}^{} \lambda_1
\ee
in probability in ${\cal
M}_{\rm R}$. This is indeed the same result as Theorem \ref{T3} in
our present setting. The
preceding
convergence is not quite sufficient to conclude:
Recalling that
$\lambda_1$ has no atoms and using Dini's theorem, we see that the
statement of the theorem
will follow if we can prove that the convergence (\ref{blocks00})
holds in the sense of
weak convergence in the space ${\cal M}_{\rm F}$. To get this
strengthening of (\ref{blocks00}), it suffices to prove the
convergence of the total masses. Therefore the proof of Theorem \ref{blocks}
will be complete once we have established the following lemma.

\begin{lemma}
\label{totalmass}
We have
$$\lim_{\varepsilon\to 0} \varepsilon\,N_{g(\varepsilon)}(]0,\infty[)
=\lambda_1(]0,\infty[)=(\Gamma(2-\gamma))^{1/(1-\gamma)},$$
in probability.
\end{lemma}

\noindent{\bf Remark.} The recent paper \cite{BBS} gives closely
related results that were
obtained independently of the present work.

\proof Write $N_{t}=N_{t}(]0,\infty[)$
to simplify notation. Then, for every $t\geq 0$ and $x\in]0,1]$, we have
$$\E[x^{N_{t}}]=\P[F_{t}(x)=1]$$
(cf formula (8) in \cite{BeLG2}). By exchangeability,
$$\P[F_{t}(x)=1]=\P[F_t(x)=F_t(1)]=\P[F_t(1-x)=0].$$
Hence, for $x\in]0,1[$,
$$\P[F_t(x)=0]=\E[(1-x)^{N_t}],$$
and it follows that
$$\P[F^\varepsilon_1(x)=0]=\E[(1-\varepsilon x)^{N_{g(\varepsilon)}}].$$
 From the convergence in distribution in Theorem \ref{smalltime},
we have for every $x>0$,
$$\limsup_{\varepsilon\to 0} \P[F^\varepsilon_1(x)=0]\leq \P[Z(1,x)=0]
=\exp-x\lambda_1(]0,\infty[).$$
We have thus obtained that, for every $x>0$,
$$\limsup_{\varepsilon\to 0}\E[(1-\varepsilon x)^{N_{g(\varepsilon)}}]
\leq \exp-x\lambda_1(]0,\infty[).$$
By standard arguments, this implies that for every $\eta>0$,
\be
\label{totalmass0}
\lim_{\varepsilon\to 0} \P[\varepsilon N_{g(\varepsilon)}<
\lambda_1(]0,\infty[) -\eta]=0.
\ee

To complete the proof, we need to verify that we have also, for every
$\eta>0$,
\be
\label{totalmasstech}
\lim_{\varepsilon\to 0}\P[\varepsilon N_{g(\varepsilon)}>
\lambda_1(]0,\infty[) +\eta]=0.
\ee
Fome now on, we fixe $\eta>0$ and we prove (\ref{totalmasstech}). We will use a
different method based on the knowledge of
the law of the process of the number of blocks in a
$\Lambda$-coalescent. For every integer $n\geq 1$, write
$N^n_t$ for the number of blocks at time $t$ in a
$\Lambda$-coalescent started initially with $n$ blocks. Then
according to Pitman \cite{Pi} (Section 3.6), the process
$(N^n_t,t\geq 0)$ is a time-homogeneous Markov chain with values
in $\{1,2,\ldots,n\}$, with only downward jumps, such that
for $2\leq k\leq b\leq n$, the
rate of jumps from $b$ to $b-k+1$ is
$$\alpha_{b,k}=\left(\begin{array}{l}b\\k\end{array}\right)
\int_{]0,1]} x^k(1-x)^{b-k}\nu(dx).$$
The total rate of jumps from $b$ is thus
$$\alpha_b=\sum_{k=2}^b \alpha_{b,k}=\int_{]0,1]} (1-(1-x)^b-b(1-x)^{b-1})
\nu(dx).$$

\begin{lemma}
\label{totaltech}
Under Assumption {\rm (A)}, we have
$$\lim_{b\to+\infty} (b^\gamma L(1/b))^{-1}\,\alpha_b= \Gamma(2-\gamma)$$
and, for every integer $k\geq 2$,
$$\lim_{b\to+\infty} (b^\gamma L(1/b))^{-1}\,\alpha_{b,k}=
{\gamma \Gamma(k-\gamma)\over k!}.$$
\end{lemma}

We leave the easy proof to the reader. Note that
\be
\label{totaltech2}
\sum_{k=2}^\infty {\gamma \Gamma(k-\gamma)\over k!}=\Gamma(2-\gamma).
\ee
This is easily proved by using the definition of the function $\Gamma$
and then an integration by parts. Similarly, we have
\be
\label{totaltech3}
\sum_{k=2}^\infty {\gamma \Gamma(k-\gamma)\over k!\Gamma(2-\gamma)}\,(k-1)=
{1\over \gamma-1}.
\ee

Let us fix $\rho\in]0,1/8[$ sufficiently small so that
$$(\Gamma(2-\gamma)^{1/(1-\gamma)}+\eta)^{1-\gamma}
<(1-6\rho)\Gamma(2-\gamma).$$
Thanks to (\ref{totaltech2}) and (\ref{totaltech3}),  we may choose
an integer $K\geq 2\vee \varepsilon_0^{-1}$
sufficiently large so that
$${1\over \Gamma(2-\gamma)}\sum_{k=2}^K {\gamma \Gamma(k-\gamma)\over
k!}\geq 1-\rho$$
and
\be
\label{totaltech4}
\sum_{k=2}^K {\gamma \Gamma(k-\gamma)\over
k!\Gamma(2-\gamma)}\,(k-1)\geq {1\over \gamma-1}-\rho.
\ee
Then, for every $k\in\{2,3,\ldots,K\}$, we may choose $\rho_k\in]0,
\gamma\Gamma(k-\gamma)/k![$ sufficiently small so that
\be
\label{totaltech5}
{1\over \Gamma(2-\gamma)}\sum_{k=2}^K (k-1)\rho_k<\rho.
\ee
Now set
$$\beta_{b,k}=\Big({\gamma \Gamma(k-\gamma)\over k!}-\rho_k\Big)
b^\gamma L({1\over b})$$
for $b\geq K$ and $k\in\{2,3,\ldots, K\}$. We also put
$$\beta_b=\sum_{k=2}^K \beta_{b,k}.$$
Notice that
\be
\label{totaltech6}
\beta_p= \sum_{k=2}^K\Big({\gamma \Gamma(k-\gamma)\over
k!}-\rho_k\Big) b^\gamma L({1\over b})\geq
(1-2\rho)\Gamma(2-\gamma)\,b^\gamma L({1\over b}).
\ee
By Lemma \ref{totaltech}, we can choose an integer $B\geq 2K$
sufficiently large so that, for every $b\geq B-K$,
$b'\in\{b,b+1,\ldots,b+K\}$ and $k\in\{2,\ldots,K\}$, one has
\be
\label{totalineq}
\beta_{b',k}\leq \alpha_{b,k}.
\ee
Denote by $(U^n_t)_{t\geq 0}$ the continuous-time Markov
chain with values in $\N$, with initial value
$U^n_0=n$, which is absorbed in the set
$\{1,\ldots,B-1\}$ and has jump rate $\beta_{b,k}$ from $b$ to $b-k+1$ when
$b\geq B$ and $k\in\{2,3,\ldots,K\}$. Fix $n\geq B$. Then thanks
to inequality (\ref{totalineq}), we can couple the Markov chains
$(U^n_t)_{t\geq 0}$ and $(N^n_t)_{t\geq 0}$ in such a way that
$$U^n_t\geq N^n_t\quad,\hbox{ for every }t\leq T^n_B:=\inf\{s:U^n_s<B\}.$$

Now it is easy to describe the behavior of the Markov chain $(U^n_t)$.
Note that for $k\in\{2,\ldots,K\}$ and $b\geq K$ the ratio
$\beta_{b,k}/\beta_b$ does not depend on $b$.
Then denote by $S_i=\xi_1+\cdots+\xi_i$ ($i=0,1,2,\ldots$) a
discrete random walk on the nonnegative
integers started from the origin and with jump distribution
$$\P[\xi_i=k-1]={\beta_{b,k}\over \beta_b}={(\gamma \Gamma(k-\gamma)/
k!)-\rho_k\over \sum_{\ell=2}^K ((\gamma \Gamma(\ell-\gamma)/
\ell!)-\rho_\ell)}\;,\qquad 2\leq k\leq K.$$
 From (\ref{totaltech4}) and (\ref{totaltech5}) we have
\be
\label{totaltech60}
\E[\xi_i]\geq {1\over \gamma-1} -2\rho.
\ee
Let ${\bf e}_0,{\bf e}_1,\ldots$ be a sequence of independent
exponential variables with mean $1$, which are also independent of
the random walk $(S_i)_{i\geq 0}$. We can construct the
Markov chain $(U^n_t)$ by setting:
$$\begin{array}{ll}
U^n_t=n&\hbox{if }0\leq t<{\displaystyle {{\bf e}_0\over \beta_n}}\\
\noalign{\smallskip}
U^n_t=n-S_1\qquad&\hbox{if }{\displaystyle {{\bf e}_0\over \beta_n}\leq
t<{\displaystyle {{\bf e}_0\over
\beta_n}+{{\bf e}_1\over \beta_{n-S_1}}}}
\end{array}
$$
and more generally,
$$U^n_t=n-S_p\qquad\hbox{if }{\displaystyle {{\bf e}_0\over
\beta_n}+{{\bf e}_1\over \beta_{n-S_1}}+\cdots
+{{\bf e}_{p-1}\over \beta_{n-S_{p-1}}}}\leq t<{\displaystyle {{\bf
e}_0\over
\beta_n}+{{\bf e}_1\over \beta_{n-S_1}}
+\cdots +{{\bf e}_{p}\over \beta_{n-S_p}}}
$$
provided $p\leq p^n_B:=\inf\{i:n-S_i<B\}$.

Recall that our goal is to prove (\ref{totalmasstech}). To this end,
note that for $a>B$,
\be
\label{totaltech7}
\P[N^n_{g(\varepsilon)}> a]\leq \P[U^n_{g(\varepsilon)}>
a]
\leq \P\Big[g(\varepsilon)\leq {{\bf
e}_0\over\beta_n}+{{\bf e}_1\over \beta_{n-S_1}}
+\cdots +{{\bf e}_{p^n_a}\over \beta_{n-S_{p^n_a}}}\Big]
\ee
where
$p^n_a:=\inf\{i:n-S_i<a\}$.

\begin{lemma}
\label{keytotal}
For $\varepsilon>0$ set
$a(\varepsilon)=(\lambda_1(]0,\infty[)+\eta)/\varepsilon$. Then,
$$\lim_{\varepsilon\to 0}\Big(\sup_{n\geq a(\varepsilon)}
\P\Big[g(\varepsilon)\leq {{\bf
e}_0\over
\beta_n}+{{\bf e}_1\over \beta_{n-S_1}}
+\cdots +{{\bf e}_{p^n_{a(\varepsilon)}}\over
\beta_{n-S_{p^n_{a(\varepsilon)}}}}\Big]\Big)=0.$$
\end{lemma}

The desired bound (\ref{totalmasstech}) immediately follows from
Lemma \ref{keytotal}. Indeed standard properties of $\Lambda$-coalescent
give
$$\P[\varepsilon N_{g(\varepsilon)}>\lambda_1(]0,\infty[)+\eta]
=\lim_{n\uparrow\infty}\uparrow \P[\varepsilon
N^n_{g(\varepsilon)}>\lambda_1(]0,\infty[)+\eta]
=\lim_{n\uparrow\infty}\uparrow \P[N^n_{g(\varepsilon)}> a(\varepsilon)]$$
and by combining (\ref{totaltech7}) and Lemma \ref{keytotal}, we see
that the latter quantity tends to $0$ as $\varepsilon\to 0$.

\noindent{\bf Proof of Lemma \ref{keytotal}:} By
(\ref{totaltech6}), we have for $a>B$,
\be
\label{totaltech8}
{{\bf
e}_0\over\beta_n}+{{\bf e}_1\over \beta_{n-S_1}}
+\cdots +{{\bf e}_{p^n_a}\over \beta_{n-S_{p^n_a}}}
\leq ((1-2\rho)\Gamma(2-\gamma))^{-1}\sum_{i=0}^{p^n_a} {{\bf e}_i\over
(n-S_i)^\gamma L({1\over n-S_i})}.
\ee
Note that
$$\E\Big[\sum_{i=0}^{p^n_a} {{\bf e}_i\over
(n-S_i)^\gamma L({1\over n-S_i})}\;\Big|\;S_i,i\geq 0\Big]
=\sum_{i=0}^{p^n_a} {1\over
(n-S_i)^\gamma L({1\over n-S_i})}.$$
Let $m\geq 2$ be an integer. For $a>B$ and $n>m a$, a trivial bound shows that
$$a^{\gamma-1}L({1\over a})\sum_{i=0}^{p^n_{ma}}{1\over
(n-S_i)^\gamma L({1\over n-S_i})}
\leq a^{\gamma-1}L({1\over a})\sum_{j=[ma]-K}^\infty {1\over j^\gamma
L({1\over j})}$$
and the right-hand side tends to $0$ as $m\to\infty$, uniformly
in $a>B$.
On the other hand, an easy argument using the law of large
numbers for the sequence $(S_i)_{i\geq 0}$ shows that, for each fixed
$m\geq 2$,
$$\lim_{a\to\infty}\Big(\sup_{n>ma}\E\Big[\Big|a^{\gamma-1}L({1\over a})
\sum_{i=p^n_{ma}}^{p^n_a}{1\over
(n-S_i)^\gamma L({1\over n-S_i})}-{1\over \E[\xi_1]}\int_1^m {dx\over
x^\gamma}\Big|\Big]\Big)=0.$$
Now recall the bound (\ref{totaltech60}) for $\E[\xi_1]$. It follows from
the preceding considerations that
\be
\label{totaltech10}
\lim_{a\to\infty} \Big(\sup_{n>a}\P\Big[a^{\gamma-1}L({1\over a})
\sum_{i=0}^{p^n_a}{1\over
(n-S_i)^\gamma L({1\over n-S_i})}> {1\over 1-3\rho}\Big]\Big)=0.
\ee

Now we can also get an estimate for the conditional variance
\ba
{\rm var}\Big(\sum_{i=0}^{p^n_a} {{\bf e}_i\over
(n-S_i)^\gamma L({1\over n-S_i})}\;\Big|\;S_i,i\geq 0\Big)
&=&\sum_{i=0}^{p^n_a} {1\over
(n-S_i)^{2\gamma} L({1\over n-S_i})^2}\\
&\leq&\sum_{j=[a-K]}^{n} {1\over j^{2\gamma} L({1\over j})^2}\\
&\leq&Ca^{1-2\gamma}L({1\over a})^{-2}
\ea
for some constant $C$ independent of $a$ and $n$. From this estimate,
(\ref{totaltech10}) and an application of the Bienaym\'e-Cebycev
inequality, we get
\be
\label{totaltech11}
\lim_{a\to\infty} \Big(\sup_{n>a}\P\Big[a^{\gamma-1}L({1\over a})
\sum_{i=0}^{p^n_a}{{\bf e}_i\over
(n-S_i)^\gamma L({1\over n-S_i})}> {1\over 1-4\rho}\Big]\Big)=0.
\ee
Recalling (\ref{totaltech8}),
we arrive at
$$\lim_{\varepsilon\to 0}
\Big(\inf_{n\geq a(\varepsilon)}
\P\Big[{{\bf
e}_0\over
\beta_n}+{{\bf e}_1\over \beta_{n-S_1}}
+\cdots +{{\bf e}_{p^n_{a(\varepsilon)}}\over
\beta_{n-S_{p^n_{a(\varepsilon)}}}} \leq
{a(\varepsilon)^{1-\gamma} L({1\over a(\varepsilon)})^{-1}\over
(1-2\rho)(1-4\rho)\Gamma(2-\gamma)}\Big]\Big)=1.$$
However, from our choice of $\rho$, we have for $\varepsilon$
sufficiently small
$$g(\varepsilon) > {a(\varepsilon)^{1-\gamma}
L({1\over a(\varepsilon)})^{-1}\over
(1-2\rho)(1-4\rho)\Gamma(2-\gamma)},$$
and this completes the proof. \QED

\noindent{\bf Remark.} It is rather unfortunate that the simple argument we
used to derive (\ref{totalmass0}) does not apply
to (\ref{totalmasstech}). On the other hand, it is interesting
to observe that the techniques involved in our proof of
(\ref{totalmasstech}) would become more complicated if we were trying to use
them to get (\ref{totalmass0}).

     \end{document}